\documentclass[3p]{elsarticle}

\usepackage{amssymb}
\usepackage{subfigure}
\usepackage{bm}
\usepackage{float}
\usepackage{algorithmic}
\usepackage{algorithm}
\usepackage{setspace}
\usepackage{booktabs}
\usepackage{enumerate}
\usepackage[misc]{ifsym}
\usepackage{url}
\usepackage{multirow}

\begin{document}

\begin{frontmatter}



\title{Molecular Sparse Representation by 3D Ellipsoid Radial Basis Function Neural Networks via $L_1$ Regularization}


\author[LSEC,SUDA]{Sheng Gui}
\address[LSEC]{LSEC, NCMIS, Academy of Mathematics and Systems Science, Chinese Academy of Sciences, Beijing 100190, China, and School of Mathematical Sciences, University of Chinese Academy of Sciences, Beijing 100049, China.}
\address[SUDA]{Department of Mathematics, Soochow University, Suzhou 215006, China.}

\author[SUDA]{Zhaodi Chen}

\author[SUDA]{Minxin Chen\corref{cor1}}
\ead{chenminxin@suda.edu.cn}
\author[LSEC]{Benzhuo Lu\corref{cor1}}
\ead{bzlu@lsec.cc.ac.cn}

\cortext[cor1]{Corresponding author}
\begin{abstract}
  In this paper, we have developed an ellipsoid radial basis function neural network (ERBFNN) and algorithm for sparse representing of a molecular shape. To evaluate a sparse representation of the molecular shape model, the Gaussian density map of molecule is approximated by ERBFNN with a relatively small number of neurons. The deep learning models were trained by optimizing a nonlinear loss function with $L_1$ regularization. Experimental results demonstrate that the original molecular shape is able to be represented with good accuracy by much fewer scale of ERBFNN by our algorithm. And our network in principle can be applied to multi-resolution sparse representation of molecular shape and coarse-grained molecular modeling.
\end{abstract}


\begin{keyword}
Molecular shape; Gaussian density map; Radial basis function neural network; Sparse representation.


\end{keyword}

\end{frontmatter}


\section{Introduction}\label{section1}
With recent advances in the field of deep learning in general, neural networks have been widely and successfully used for different tasks in computer vision, including object detection \cite{2013arXivSermanet,Girshick2015FastRCNN}, image classification \cite{He2016DeepResidual,Szegedy2016arXivInception} and semantic segmentation \cite{Ronneberger2015UNet,2012arXivEfficientInference}. As a special class of feedforward neural networks (FNNs), radial basis function (RBF) networks have certain advantages over other types of FNNs, such as simpler network structures and faster training process. Due to good approximation capabilities, single-output RBF networks are usually utilized to model nonlinear functions in engineering applications. In practice, learning of RBF networks includes two assignments that determining the RBF network structure and optimizing the adaptable parameters (such as, centers and their radii of RBF neurons, and linear output weights). The \cite{Peng2007RBF} has described, if the two tasks are meanwhile taken into consideration, it becomes a mixed integer programming of hard problem. Owing to the lack of promising methods to address this integrated problem, the two tasks are solved respectively in many learning algorithms of RBF networks \cite{Schwenker2001RBF}. In this case, the network structure is determined in advance, and the parameters are then trained by algorithms of supervised learning. As well known that directly optimizing the empirical risk may lead to an overfitting problem, which causes poor generalization capability. To tackle this issue, the regularization techniques, such as $L_1$ and $L_2$ \cite{Andrew07scalabletraining,scholkopf2007L1,Perkins2003Grafting,shen2013FullyLoss,Vidaurre2010L1} regularization, are wealthy in modern machine learning. Fundamentally, a regularization term is added to the empirical risk to penalize over-complicated solutions. $L_1$ regularization is implemented by appending a weighted $L_1$ norm of the parameter vector to the loss function, which ensures the sum of the absolute values of the parameters to be small. While $L_2$ regularization uses the $L_2$ norm, which encourages the sum of the squares of the parameters to be small. There has been an increasing interest in $L_1$ regularization because of its advantages over $L_2$ regularization \cite{Qian2017SpareRBF}. For example, $L_1$ regularization usually produces sparse parameter vectors in which many parameters are closed to zero. Thus, the more sparse solution can be obtained. In particular, if one could use deep learning \cite{Girosi1998equivalence} to representing sparsely the shape for biomolecules, that would result in faster and more efficient way for such as molecular docking, alignment, drug design and multiscale modeling.

Biomolecules such as proteins are the fundamental functional units of life activities. Geometric modeling of biomolecules plays an important role in the fields of computer-aided drug design and computational biology \cite{Liao2015Atom}. In the computer-aided drug design field, biomolecular shape has been a vital issue of considerable interest for many years, for instance, shape-based docking problems \cite{mcgann2003gaussian}, molecular shape comparisons \cite{grant1996a}, calculating molecular surface areas \cite{Liu2015Parameterization,Weiser1999Gaussian}, coarse-grained molecular dynamics \cite{Rudd2000coarse}, and the generalized Born models\cite{yu2006what}, etc. And biomolecules geometric shape (especially molecular surface) is prerequisite for using boundary element method (BEM) and finite element method (FEM) in the implicit solvent models \cite{Lu2008Recent}. Considering the highly complex and irregular shape of a molecule, new challenges arise in simulations involving extremely large biomolecular \cite{Chen2013Advances} (e.g., viruses, biomolecular complexes etc.). And the efficient representation of the molecular shape (as well as the "molecular surface" or "molecular volume") for large real biomolecule with high quality remains a critical topic \cite{Lu2008Recent}.

The molecular shape is defined in various senses \cite{CONNOLLY1983Analytical,Gerstein2001Protein,DUNCAN1993ShapeAnalysis}. For molecular volume, the Gaussian density map is a suitable representation of the molecular shape, since the Gaussian density maps provide a realistic representation of the volumetric synthetic electron density maps for the biomolecules \cite{DUNCAN1993ShapeAnalysis}. For molecular surface, there are four important biomolecular surfaces: van der Waals (VDW) surface, solvent accessible surface (SAS) \cite{lee1971the}, solvent excluded surface (SES) \cite{RICHARDS1977AREAS} and Gaussian surface. The van der Waals surface is the smallest envelope enclosing a collection of spheres representing all the atoms in the system with their van der Waals radii. The SAS  \cite{LEE1971INTERPRETATION}  is the trace of the centers of probe spheres rolling over the van der Waals surface. The SES \cite{RICHARDS1977AREAS} is the surface traced by the inward-facing surface of the probe sphere. The Gaussian surface \cite{Weiser1999Gaussian,Zhang2006Qualitymeshing} is a level-set of the summation of the spherically symmetrical volumetric Gaussian density distribution centered at each atom of the biomolecular system. In 2015, Liu et al. presented that the VDW surface, SAS, and SES can be approximated well by the Gaussian surface with proper parameter selection in sense of parameterizations for molecular Gaussian surface \cite{Liu2015Parameterization}. Comparing with VDW surface, SAS and SES, the Gaussian surface  is smooth and has been widely used in many problems in computational biology \cite{mcgann2003gaussian,grant1996a,Weiser1999Gaussian,Rudd2000coarse,yu2006what}.  Thus, in this paper, we adopt ellipsoid RBF neural network to approximate to the Gaussian density maps of molecular shape. The Gaussian density maps and the Gaussian surface descriptions of the specific forms will be given in the next section.

For Gaussian density maps, the volume Gaussian function is constructed by a summation of Gaussian kernel functions, whose number depends on the total number of atoms in the molecular. Thus, the computational cost for biomolecular surface construction increases as the atom number (number of Gaussian kernel functions) becomes progressively larger. It leads to a significant challenge for their analysis and recognition. In case of large biomolecules, the number of kernels in their definition of Gaussian molecular surface may achieve millions. In 2015, Zhang et al. \cite{Liao2015Atom} put forward an atom simplification method for the biomolecular structure based on Gaussian molecular surface. This method contains two main steps. The first step eliminates the low-contributing atoms. The second step optimizes the center location, the radius and the decay rate of the remaining atoms based on gradient flow method.

In the area of computer aided geometric design, the Gaussian surface is a classical implicit surface representing method. Over that last few decades, there are a mount of works focusing on the implicit surface reconstruction problem, and various approaches have been presented. J.C. Carr \cite{Carr2001Reconstruction} proposed a method to reconstruct an implicit surface with RBFs and performed a greedy algorithm to append centers with large residuals to decrease the number of basis functions. However, the result of this method is not sparse enough. M. Samozino \cite{Samozino2006RVC} presented a strategy to put the RBFs’ centers on the Voronoi vertices. This strategy, firstly picks a user-specified number of centers by filtering and clustering from a subset of Voronoi vertices, then gets the reconstructed surface by solving a least-square problem. However, it leads to larger approximation error on the surface while approximating the surface and center points equally. In 2016, Chen \cite{2016LiSparse} et al. proposed a model of sparse RBF surface representations. They constructed the implicit surface based on sparse optimization with RBF. And the initial Gaussian RBF is on the medial axis of the input model. They have solved the RBF surface by sparse optimization technique. Sparse optimization has become a very popular technique in many active fields, for instance, signal processing and computer vision, etc \cite{Elad2010}. This technique has been applied in linear regression \cite{Atkinson1992Subset}, deconvolution \cite{TAYLOR1979DECONVOLUTION}, signal modeling \cite{RISSANEN1978MODELING}, preconditioning \cite{Grote1997Parallel}, machine learning \cite{Girosi1998equivalence}, denoising \cite{Chen1998Atomic}, and regularization \cite{Daubechies2004sparsity}. In the last few years, sparse optimization also has been applied in geometric modeling and graphics problems (refer to a review \cite{Xu2015}).

In this paper, based on the structure of RBF network, we propose an ellipsoid RBF neural network for reducing the number of kernels in the definition of Gaussian surface while preserving the shape of the molecular surface. We highlight several differences and main contributions between our method and previous $L_1$ optimization methods with shape representation:
\begin{enumerate}
    \item Compared with other works, our focus is mainly on reducing the number of kernels in  Gaussian density maps by pruning useless ellipsoid RBF neuron through $L_1$ regularization;
    \item The loss function of our model is a complicated nonlinear function with respect to the locations, sizes, shapes and orientations of RBFs;
    \item Different initializations and training network algorithms are proposed for solving the corresponding optimization problem in our model.
\end{enumerate}

The remainder of this paper is organized as follows. Section \ref{section2} reviews some preliminary knowledge about volumetric electron density maps, Gaussian surface, ellipsoid Gaussian RBF and ellipsoid RBF network, then presents our model together with an algorithm for representing the Gaussian density maps sparsely. The experimental results and comparisons are demonstrated in section \ref{section3}. We conclude the paper in section \ref{section4}.

\section{Methods}\label{section2}
\subsection{Brief review of  volumetric electron density maps, Gaussian surface, ellipsoid Gaussian RBF and ellipsoid RBF network}
\subsubsection{Volumetric electron density maps}
Volumetric electron density maps are often modelled as the volumetric Gaussian density maps $\phi : \mathbb{R}^3 \rightarrow \mathbb{R}$. The definition of the volumetric Gaussian density maps is as follows,
 \begin{equation}
\label{GaussKernel}
\phi ( \mathbf{x}) = \sum \limits _{i = 1}^N e^{-d(\| \mathbf{x} -  \mathbf{x}_i \|^2 - r_i^2)},
\end{equation}
where the parameter $d$ is positive and controls the decay rate of the kernel functions, $\mathbf{x}_i$ and $r_i$ are the location and radius of atom $i$.
\subsubsection{Gaussian surface}
\begin{figure}[H]
  \centering
  \subfigure[]{
  \label{Gaussian:PQR}
  \includegraphics[scale = 0.24]{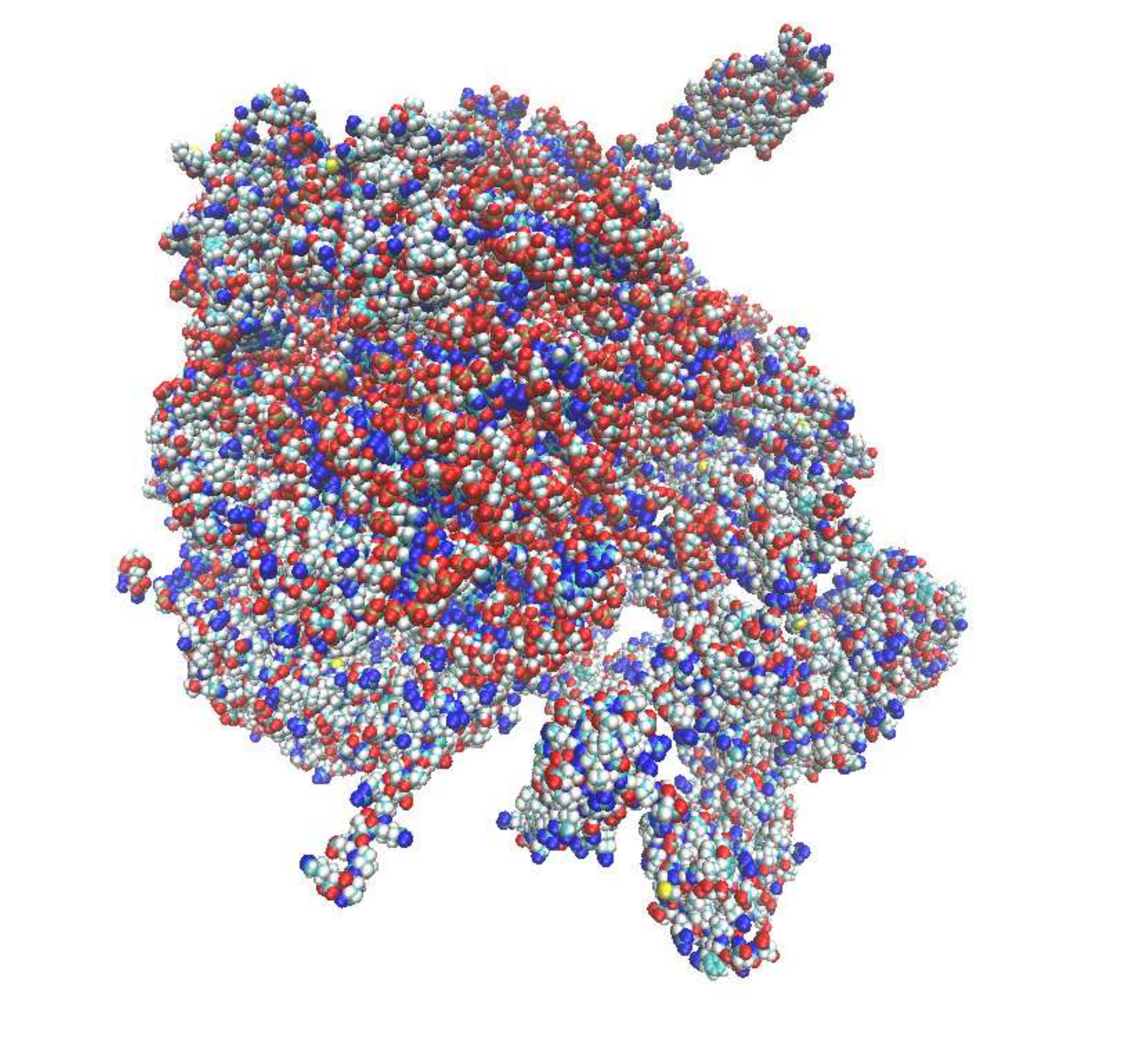}
  }
  \subfigure[]{
  \label{Gaussian:Surface}
  \includegraphics[scale = 0.22]{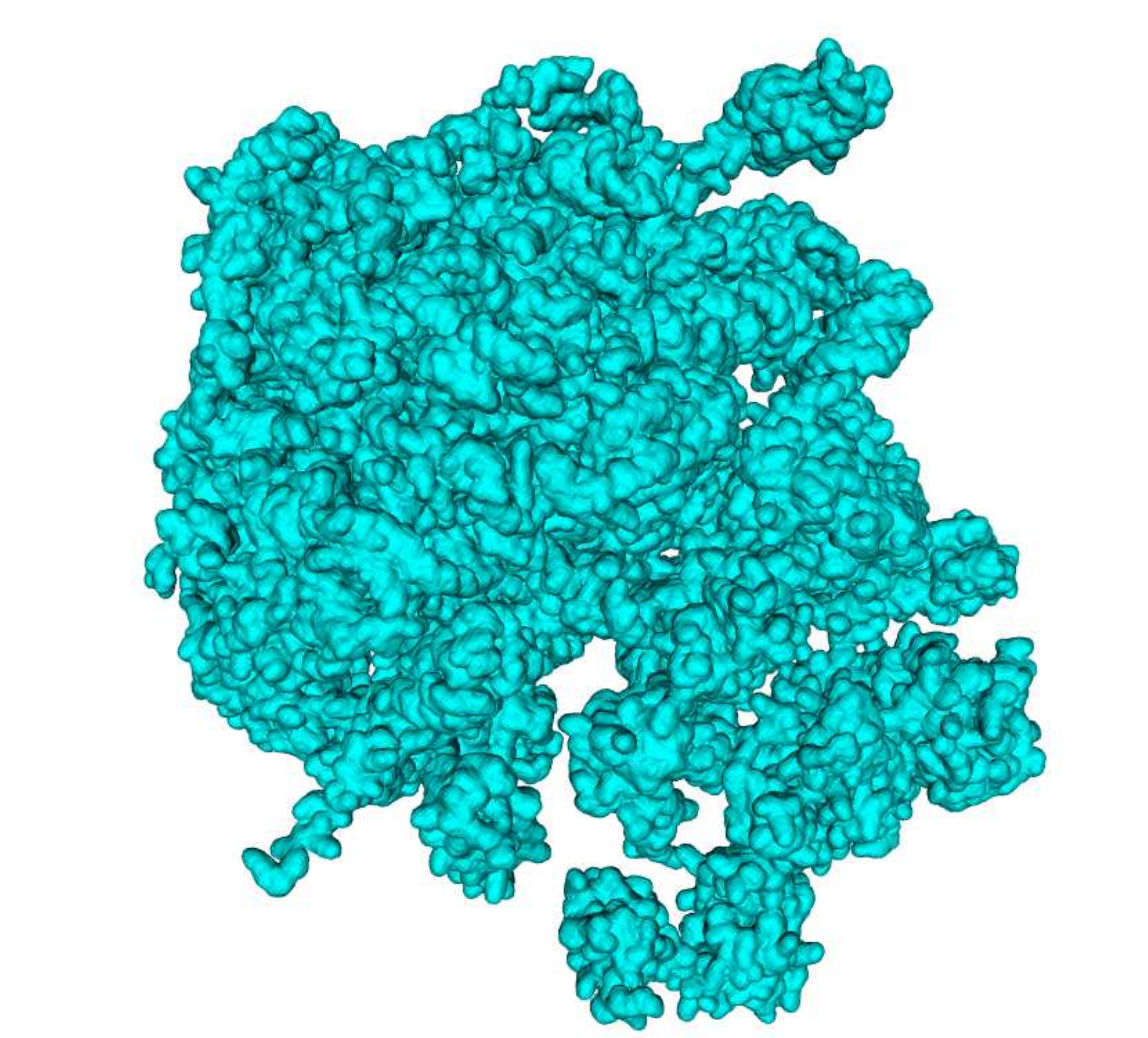}
  }

  \caption{An example of Gaussian molecular surface via VCMM\cite{Bai2014VCMM}. (a) shows the VDW surface, and (b) shows the Gaussian molecular surface generated by TMSmesh \cite{Chen2012Triangulated,Chen2011TMSmesh,Liu2018EFFICIENT} with parameter $d$ and $c$ being setting as 0.9 and 1.0, respectively. All coordinates and corresponding radii are drawn from the PQR file that is transformed from the PDB file, using the PDB2PQR tool \cite{Dolinsky2004PDB2PQR}.}
  \label{GaussianSurface}
\end{figure}

The Gaussian surface is defined as a level set from volumetric synthetic electron density maps,
\begin{equation}
\label{GaussLevelSet}
\left\{ \mathbf{x}  \in \mathbb{R}^3,\phi \left( { \mathbf{x}} \right) = c \right\},
\end{equation}
where  $c$ is the isovalue, and it controls the volume enclosed by the Gaussian density maps. Fig. \ref{GaussianSurface} shows an example of a Gaussian surface. This molecule contains the entire 70S ribosome, including the 30S subunit (16S rRNA and small subunit proteins), 50S subunit (23S rRNA, 5S rRNA, and large subunit proteins), P- and E-site tRNA, and messenger RNA. This molecule is obtained from 70S ribosome3.7A model140.pdb.gz on \url{http://rna.ucsc.edu/rnacenter/ribosome downloads.html}. Fig. \ref{Gaussian:PQR} shows all the atoms in the molecule, and Fig. \ref{Gaussian:Surface} shows the corresponding Gaussian surface.

\subsubsection{Ellipsoid Gaussian RBF}
The RBF is written as $\xi_i (\mathbf{x}) = \xi(\|\mathbf{x} - \mathbf{c}_i\|)$, where $\xi(\mathbf{x})$ is a nonnegative function defined on $[0,\infty)$, the $\mathbf{c}_i$ is center location of $i$th basis function. The RBF has basic properties as following, $\xi(0) = 1$  and $lim _{\mathbf{x} \rightarrow + \infty} \xi (\mathbf{x}) = 0$.  A typical choice of RBF is Gaussian function
\begin{equation}\label{GaussianRBF}
  \xi (\mathbf{x}) = e^{-\mathbf{x}^2}.
\end{equation}
In addition, there are other RBF including thin plate spline RBF, e.g., $\xi (r) = r^2ln(r)$ for $r \in \mathbb{R}$ .

Compared with other RBF, we put forward ellipsoid RBF with parameters that respect to the locations, sizes, shapes and orientations. The ellipsoid Gaussian RBF can be rewritten as,
\begin{equation}\label{ERBF}
  \psi(\mathbf{x}) = e^{-\|\mathbf{D}^{\frac{1}{2}}\bm{\Theta}(\alpha, \beta, \gamma)\left(\mathbf{x} - \mathbf{c}\right)\|^2},
\end{equation}
where $\mathbf{c} = (c_1, c_2, c_3)^T \in \mathbb{R}^3$ is the center of the ellipsoid Gaussian RBF, $\mathbf{D} =  diag(d_1,d_2,d_3)$, where $d_i, i = 1,2,3$ defines the length of ellipsoid along three main axis,
$\bm{\Theta}(\alpha, \beta, \gamma)$ is the total rotation matrix, and it is equal to the product of rotation matrices from three directions
\begin{equation}
\label{RotataionMatirx}
\bm{\Theta}(\alpha, \beta, \gamma) = \bm{\Theta}_z(\gamma ) \cdot \bm{\Theta}_y(\beta )\cdot \bm{\Theta}_x(\alpha ),
\end{equation}
and $\bm{\Theta}_x(\alpha )$ is a rotation matrix of x direction:
\begin{equation}
\bm{\Theta_x}\left(\alpha \right)=\left[
 \begin{array}{ccc}
    {1} & {0} & {0} \\
    {0} & {\cos \alpha} & {-\sin \alpha} \\
    {0} & {\sin \alpha} & {\cos \alpha}
 \end{array}
\right],
\end{equation}
$\bm{\Theta}_y(\beta )$ is a rotation matrix of y direction:
\begin{equation}
\bm{\Theta}_y\left(\beta\right)=\left[
 \begin{array}{ccc}
    {\cos \beta} & {0} & {\sin \beta} \\
    {0} & {1} & {0} \\
    {-\sin \beta} & {0} & {\cos \beta}
 \end{array}\right],
\end{equation}
$\bm{\Theta}_z(\gamma )$ is a rotation matrix of z direction:
\begin{equation}
\bm{\Theta}_z\left(\gamma\right)=\left[
 \begin{array}{ccc}
    {\cos \gamma} & {-\sin \gamma} & {0} \\
    {\sin \gamma} & {\cos \gamma} & {0} \\
    {0} & {0} & {1}
 \end{array}\right],
\end{equation}
so that $\bm{\Theta}(\alpha, \beta, \gamma)$ is equal to:
\begin{equation}
\left[
\begin{array}{ccc}
{\cos \beta \cos \gamma} & {-\cos \alpha \sin \gamma +  \sin \alpha \sin \beta \cos \gamma} & {\sin \alpha \sin \gamma + \cos \alpha \cos \gamma \sin \beta} \\
{\cos \beta \sin \gamma} & {\cos \alpha \cos \gamma + \sin \alpha \sin \beta \sin \gamma} & {-  \sin \alpha \cos \gamma + \cos \alpha \sin \beta \sin \gamma} \\
{-\sin \beta } & {\cos \beta \sin \alpha} & {\cos \alpha \cos \beta}
\end{array}
\right].
\end{equation}

\subsubsection{Ellipsoid RBF Networks}
The RBF network is a special FNN consisting of three layers:
\begin{itemize}
	\item an input layer
	\item a hidden layer with a nonlinear activation function
	\item a linear output layer
\end{itemize}

The choice of activation function is the ellipsoid Gaussian function. For an input $\mathbf{x} \in \mathbb{R}^3$, the output of the ellipsoid RBF network is calculated by
\begin{equation}
\Psi(\mathbf{x}) =  \sum \limits _{i = 1}^{N} w_i \psi_i(\mathbf{x}) = \sum \limits _{i = 1}^{N} w_i e^{-\|\mathbf{D}_i^{\frac{1}{2}}\bm{\Theta}_i(\alpha_i, \beta_i, \gamma_i)\left(\mathbf{x} - \mathbf{c}_i\right)\|^2},
\end{equation}
where $\mathbf{c}_i = \left[c_{i1}, c_{i2}, c_{i3}\right]^{\top} \in \mathbb{R}^3$ is the $i$th ellipsoid RBF center of hidden layer, $\mathbf{D}_i = diag(d_{i1},d_{i2},d_{i3})$ represents the lengths of corresponding ellipsoid RBF along three main axes of hidden layer, $\bm{\Theta}_i(\alpha_i, \beta_i, \gamma_i)$ is a rotation matrix of the $i$th neuron. $w_i$ is the output weight between the $i$th hidden neuron and the output node. And $\|\cdot\|$ is the $L_2$ norm of vector.

Denote the parameters (i.e., the weights connecting the neuron to the output layer, lengths of centers, centers coordinate and rotation angles) of the hidden neuron by $\bm{\sigma} = \left[\mathbf{w}, \mathbf{d}, \mathbf{c}, \bm{\alpha}, \bm{\beta}, \bm{\gamma} \right]^\top \in \mathbb{R}^{10N}$. And the descriptions of the specific forms of $\bm{\sigma}$ will be given in the following section. Assume the training data set is given by $\left\{(\mathbf{x}_m, \mathbf{y}_m)| \mathbf{x}_m \in \mathbb{R}^3, \mathbf{y}_m \in \mathbb{R}, m = 1,2, \cdots, M\right\}$, where $\mathbf{x}_m$ is the $m$th input pattern and $\mathbf{y}_m$ is the desired output value for the $m$th input pattern. The actual output vector can be calculated by

\begin{equation}\label{outputVector}
  \mathbf{\hat{Y}} = \mathbf{H} \mathbf{w}
\end{equation}
where
\begin{equation}       
\mathbf{H} = \left[                 
  \begin{array}{cccc}   
    e^{-\|\mathbf{D}_1^{\frac{1}{2}}\bm{\Theta}_1(\mathbf{x}_{1} - \mathbf{c}_1)\|_2^2} & e^{-\|\mathbf{D}_2^{\frac{1}{2}}\bm{\Theta}_2(\mathbf{x}_{1} - \mathbf{c}_2)\|_2^2}&\cdots & e^{-\|\mathbf{D}_N^{\frac{1}{2}}\bm{\Theta}_N(\mathbf{x}_{1} - \mathbf{c}_N)\|_2^2}\\  
     e^{-\|\mathbf{D}_1^{\frac{1}{2}}\bm{\Theta}_1(\mathbf{x}_{2} - \mathbf{c}_1)\|_2^2} & e^{-\|\mathbf{D}_2^{\frac{1}{2}}\bm{\Theta}_2(\mathbf{x}_{2} - \mathbf{c}_2)\|_2^2}&\cdots & e^{-\|\mathbf{D}_N^{\frac{1}{2}}\bm{\Theta}_N(\mathbf{x}_{2} - \mathbf{c}_N)\|_2^2}\\  
    \vdots &  \vdots  &  & \vdots \\
    e^{-\|\mathbf{D}_1^{\frac{1}{2}}\bm{\Theta}_1(\mathbf{x}_{M} - \mathbf{c}_1)\|_2^2}& e^{-\|\mathbf{D}_2^{\frac{1}{2}}\bm{\Theta}_2(\mathbf{x}_{M} - \mathbf{c}_2)\|_2^2} & \cdots & e^{-\|\mathbf{D}_N^{\frac{1}{2}}\bm{\Theta}_N(\mathbf{x}_{M} - \mathbf{c}_N)\|_2^2}\\  
  \end{array}
\right]_{M \times N},                
\end{equation}

$\mathbf{\hat{Y}} = \left[\mathbf{\hat{y}}_1, \mathbf{\hat{y}}_2, \cdots, \mathbf{\hat{y}}_M\right]^{\top}$ is an output value vector for $M$ input patterns, $\mathbf{w} = \left[w_1, w_2, \cdots, w_N\right]^{\top}$ is a $N$ vector, $w_k$ is the weight connecting the $k$th hidden neuron to the output layer. The error vector is defined as
\begin{equation}\label{errorVector}
  \mathbf{e} = \left[e_1, e_2, \cdots, e_N\right]^{\top}
\end{equation}
with $e_i = \hat{y}_i - y_i$.
\begin{figure}[H]
  \centering
  \includegraphics[scale = 0.4]{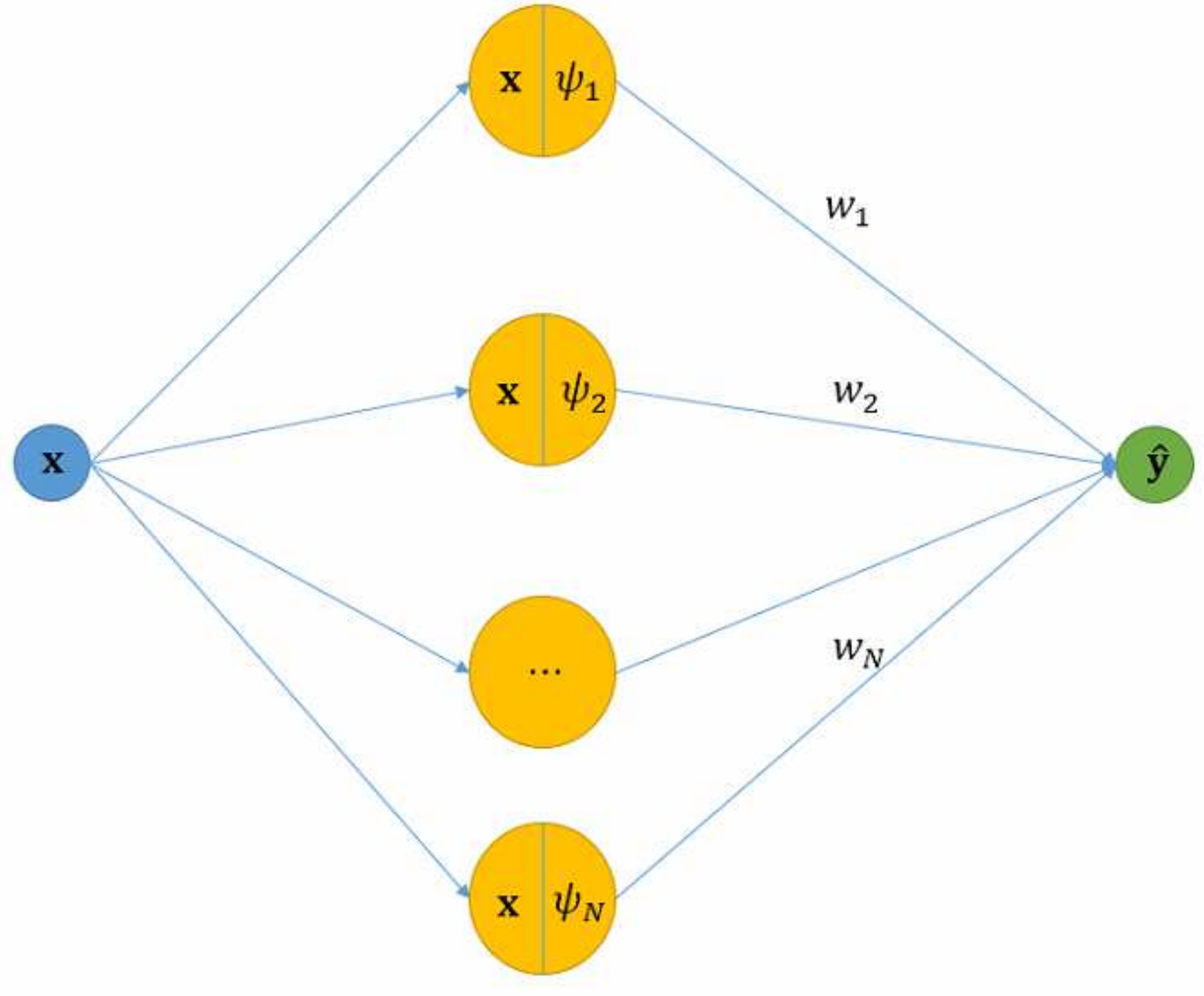}\\
    \caption{The structure of ellipsoid RBF neural network.}\label{ERBFnetwork}
\end{figure}

\subsection{Model and algorithm}
\subsubsection{Modeling with ellipsoid RBF network}
The major goal of the study was to create the sparse representation of Gaussian molecular model by the ellipsoid RBF neural network. According to the definition of volumetric electron density maps and structure of the ellipsoid RBF network, the loss function of representing sparsely Gaussian molecular is as follows,
\begin{equation}\label{MODELF}
\mathcal{L}(\bm{\sigma}) = \rho_1 \cdot \left[ \|\mathbf{w}\|_1 +  \|\mathbf{d}\|_1 \right] + \rho _2 \cdot E_s(\bm{\sigma}),
\end{equation}
and corresponding constrained condition is
\begin{equation}\label{constrainedcondition}
\mathbf{w}  \ge 0,  \mathbf{d}  \ge 0, \mathbf{w} \in \mathbb{R}^{N}, \mathbf{d} \in \mathbb{R}^{3N}.
\end{equation}

The first term in Eq. \ref{MODELF} is a $L_1$ regularization term to reduce term both network complexity and overfitting. The formulate is as follows,
\begin{equation}\label{L1}
 \|\mathbf{w}\|_1 + \|\mathbf{d}\|_1 = \sum \limits _{i = 1}^{N} |w_i| + \sum \limits _{i = 1}^{N} \sum \limits _{j = 1}^{3} |d_{ij}| ,
\end{equation}
where $\mathbf{w}=\left[w_{1,} w_{2}, \cdots, w_{N}\right]$, $\mathbf{d}=\left[d_{11}, d_{12}, d_{13}, \cdots, d_{N1}, d_{N2}, d_{N3}\right]$.

The second $E_s(\bm{\sigma})$ is density error of the sparsely representing molecule and original molecule at training points set $\mathbf{x}_m, m = 1, 2, \cdots, M$.  We have
\begin{equation}
E_s(\bm{\sigma}) = \sum \limits _{m=1}^{M}\left[\Psi(\mathbf{x}_m;\bm{\sigma}) - \phi(\mathbf{x}_m)\right]^{2}  = \sum \limits _{m=1}^{M}\left[ \sum \limits _{i = 1}^{N} w_i  e^{-\|\mathbf{D}_i^{\frac{1}{2}}\bm{\Theta}_i(\mathbf{x}_m - \mathbf{c}_i)\|_2^2} - \phi(\mathbf{x}_m)\right]^{2},
\end{equation}
where $\Psi(\mathbf{x})$ is an ellipsoid RBF neuron network. $\phi(\mathbf{x})$ is the volumetric electron density map (Eq. \ref{GaussKernel}).
It is to be approximated by $\Psi(\mathbf{x})$. $\mathbf{x}_m = (x_{m1}, x_{m2}, x_{m3})^\top \in \mathbb{R}^3$ is the $m$th training point. $\mathbf{c}_i = (c_{i1}, c_{i2}, c_{i3})^\top\in \mathbb{R}^3$ is the center of the $i$th activation function of the ellipsoid RBF. $\mathbf{D}_i = diag(d_{i1},d_{i2},d_{i3})$ define the lengths of ellipsoid along three main axis. $\bm{\Theta}_i$ is a rotation matrix. The $\alpha_i,\beta_i,\gamma_i$ are rotation angles of the $i$th activation function of the ellipsoid RBF neuron, $i = 1,2,\cdots, N$, $m = 1,2,\cdots,M$. $N$ is the number of the ellipsoid RBF neurons. $M$ is the number of the training point set. The $\bm{\sigma}$ is the network parameter, the formula is
\begin{equation}\label{X}
\bm{\sigma} = \left[\mathbf{w}, \mathbf{d}, \mathbf{c}, \bm{\alpha}, \bm{\beta}, \bm{\gamma} \right]^\top,
\end{equation}
where $\mathbf{w}=\left[w_{1,} w_{2}, \cdots, w_{N}\right]$, %
$\mathbf{d}=\left[d_{11}, d_{12}, d_{13}, \cdots, d_{N1}, d_{N2}, d_{N3}\right]$,
$\mathbf{c}=\left[\mathbf{c}_{1}, \mathbf{c}_{2}, \cdots,\mathbf{c}_{N}\right]$, $\bm{\alpha}=\left[\alpha_{1,} \alpha_{2}, \cdots, \alpha_{N}\right]$, $\bm{\beta}=\left[\beta_{1,} \beta_{2}, \cdots, \beta_{N}\right]$, $\bm{\gamma}=\left[\gamma_{1,} \gamma_{2}, \cdots, \gamma_{N}\right]$.


The two parameters $\rho _1 > 0$ and $\rho _2> 0$ are used to balance the two targets: accuracy ($E_s$) and sparsity ($L_1$-regularization). And the constrained conditions are explained as follows, (i) $\mathbf{w}>0$ indicates that the corresponding ellipsoid Gaussian RBF is nonnegative which means each RBF in $\psi$ can be seen as a new real physical atom with ellipsoid shape. (ii) $\mathbf{d} \geq 0$ implies  the activation function is zero at infinity, which is consistent with the fitted function $\phi$. In order to transform the Eq. \ref{MODELF} and Eq. \ref{constrainedcondition}  to an unconstrained loss function, we do the following substitution,
\begin{equation}
\left\{ \begin{array}{ll}
w_i  = \tilde{w_i}^2 , \\
d_{iq}  = \tilde{d}_{iq}^2, \\
i = 1,2,\cdots,N, q = 1,2,3,
\end{array} \right.
\end{equation}
and corresponding $\mathbf{\tilde{D}}_i = diag(\tilde{d}_{i1}^2,\tilde{d}_{i2}^2,\tilde{d}_{i3}^2)$. For simplicity, we still use $w_i,d_{qi},\mathbf{D}_i$ to denote $\tilde{w}_i,\tilde{d}_{iq},\mathbf{\tilde{D}}_i$.

Thus, the loss function of the ellipsoid RBF network for representing sparsely a molecule is:
\begin{equation}\label{Model2}
\mathcal{L}(\bm{\sigma}) = \rho_1 \cdot \left[\|\mathbf{w}\|_1 + \|\mathbf{d}\|_1 \right] + {\rho_2} \cdot \sum_{m=1}^{M}\left[\sum_{i=1}^{N} w_i^2 \cdot e^{-\| \mathbf{D}_i\bm{\Theta}_i(\mathbf{x}_{m} - \mathbf{c}_i)\|^2_2} - \phi\left(\mathbf{x}_{m}\right)\right]^{2}.\\
\end{equation}

\subsubsection{Overview}
\begin{figure}[H]
  \centering
  \includegraphics[scale = 0.35]{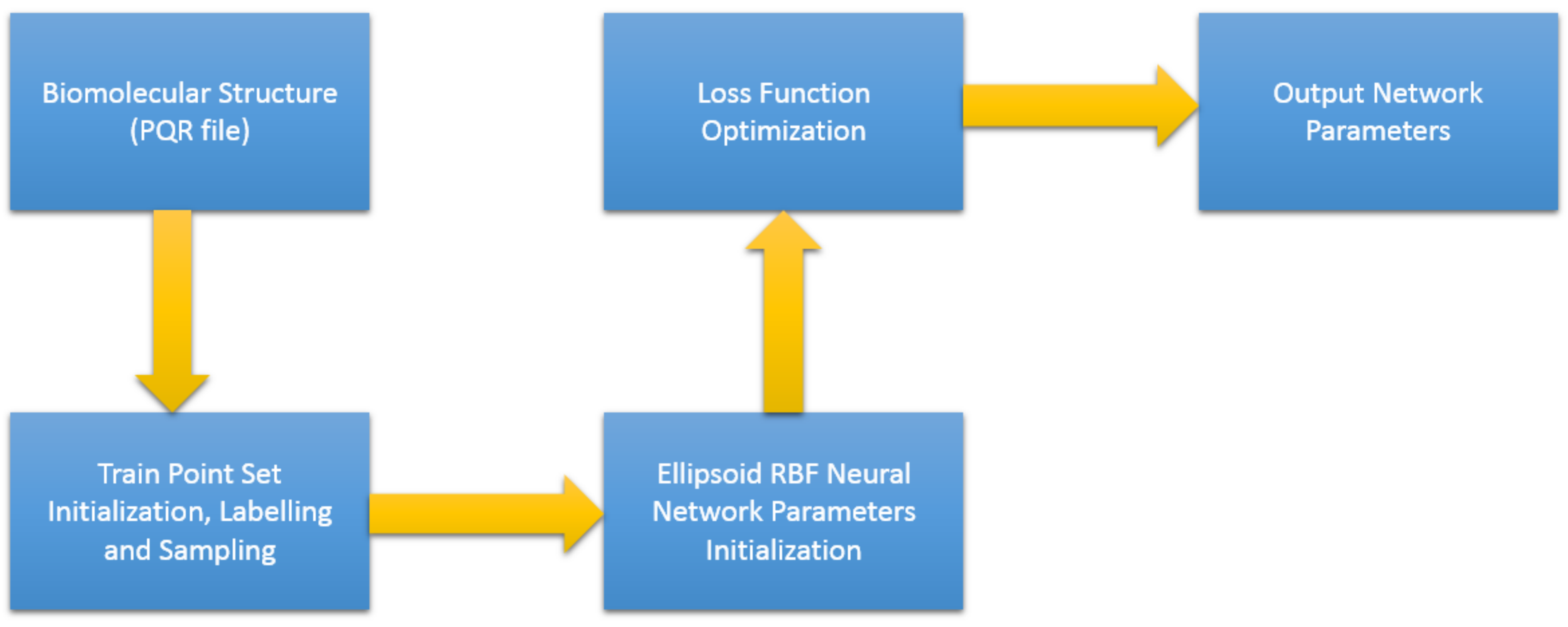}\\
  \caption{The pipline of our algorithm and results in each step.}\label{OverView}
\end{figure}

In this section, we describe the algorithms to construct sparse representation of Gaussian molecular density by the ellipsoid RBF neural network. The inputs of our method are PQR files which include a list of centers and radii of atoms. The output of our method are network parameters which contain the centers, the lengths, the rotation angles of the ellipsoid RBF neural network and the weights connecting the hidden neurons to the output layer. The algorithm outline is as follows: first,  set the training points $\mathbf{x}_m, m=1,...,M$ and label corresponding value $\phi(\mathbf{x}_m)$ . Second, initialize the ellipsoid RBF network (i.e., the number of neuron, the parameters of the ellipsoid RBF neural network). Third, optimize the loss function in Eq. \ref{Model2} using an ADAM algorithm \cite{2014arXivADAM} to minimizing the sparsity and error terms in Eq. \ref{Model2} alternatively. 
Fig. \ref{OverView} demonstrates the process of our algorithm. The result shows that, using our method, the original Gaussian surface is approximated well by a summation of much fewer ellipsoid Gaussian RBFs.

\subsubsection{Training points set initialization and labelling}
\label{Yinitialization}
\begin{figure}[H]
  \centering
  \subfigure[]{
  \label{point:surface}
  \includegraphics[scale = 0.35]{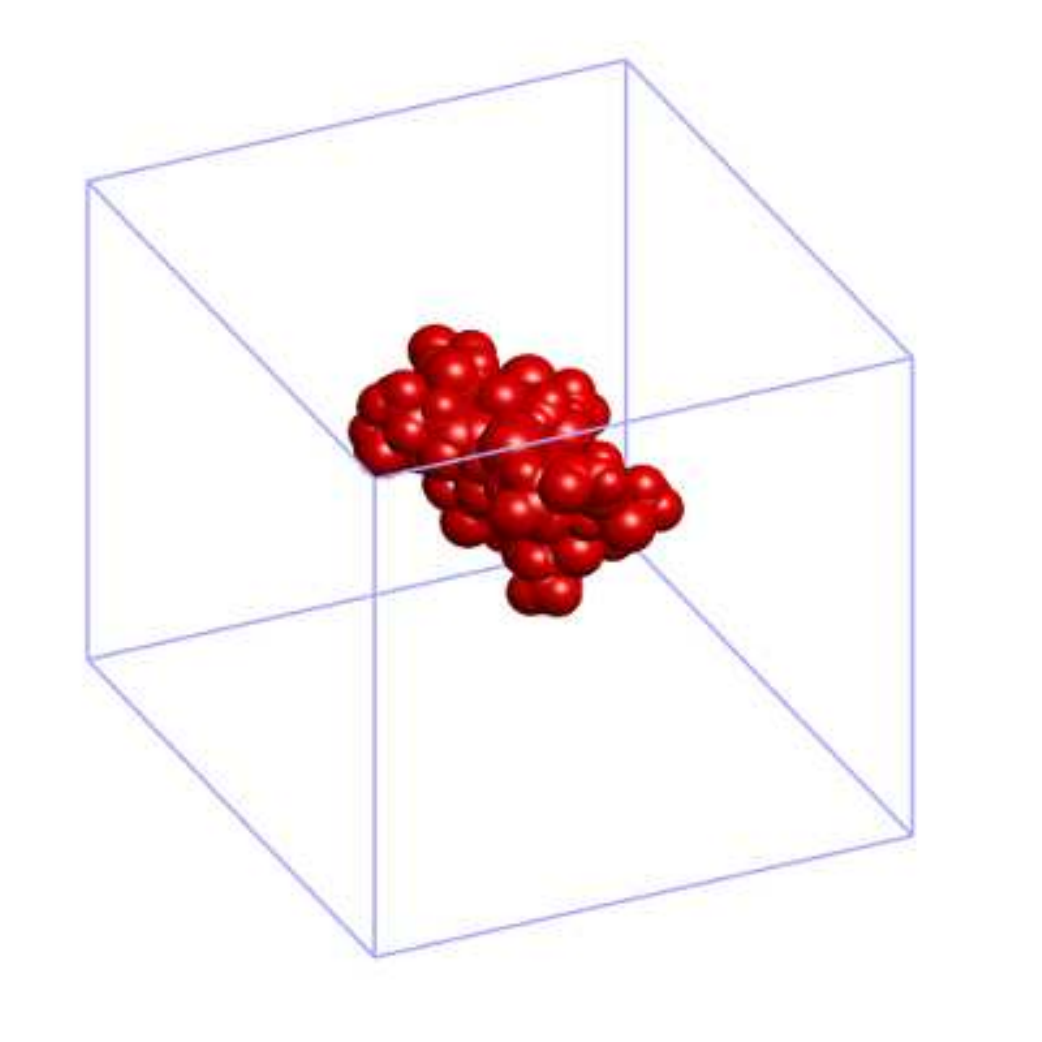}}
  \subfigure[]{
  \label{point:box}
  \includegraphics[scale = 0.35]{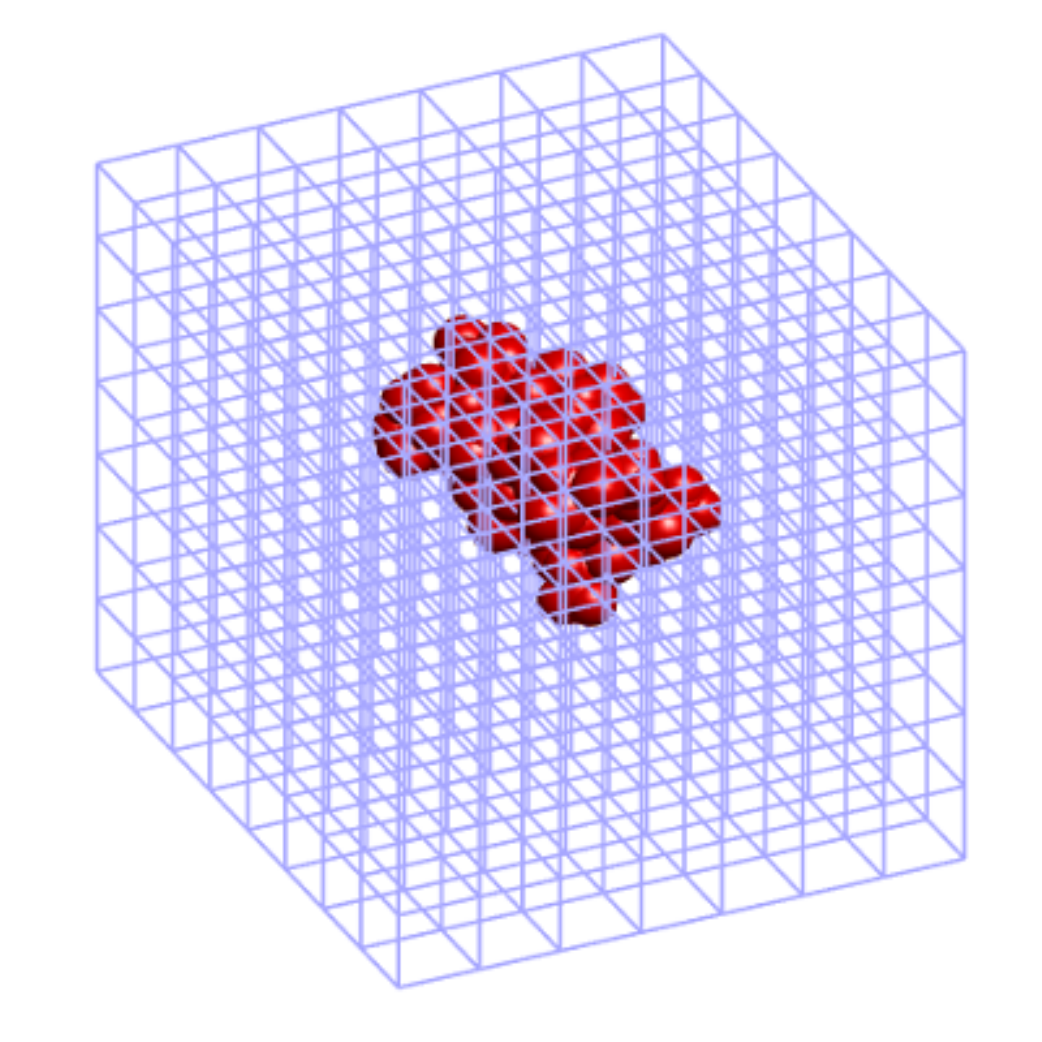}}
  \subfigure[]{
  \label{point:constraint}
  \includegraphics[scale = 0.35]{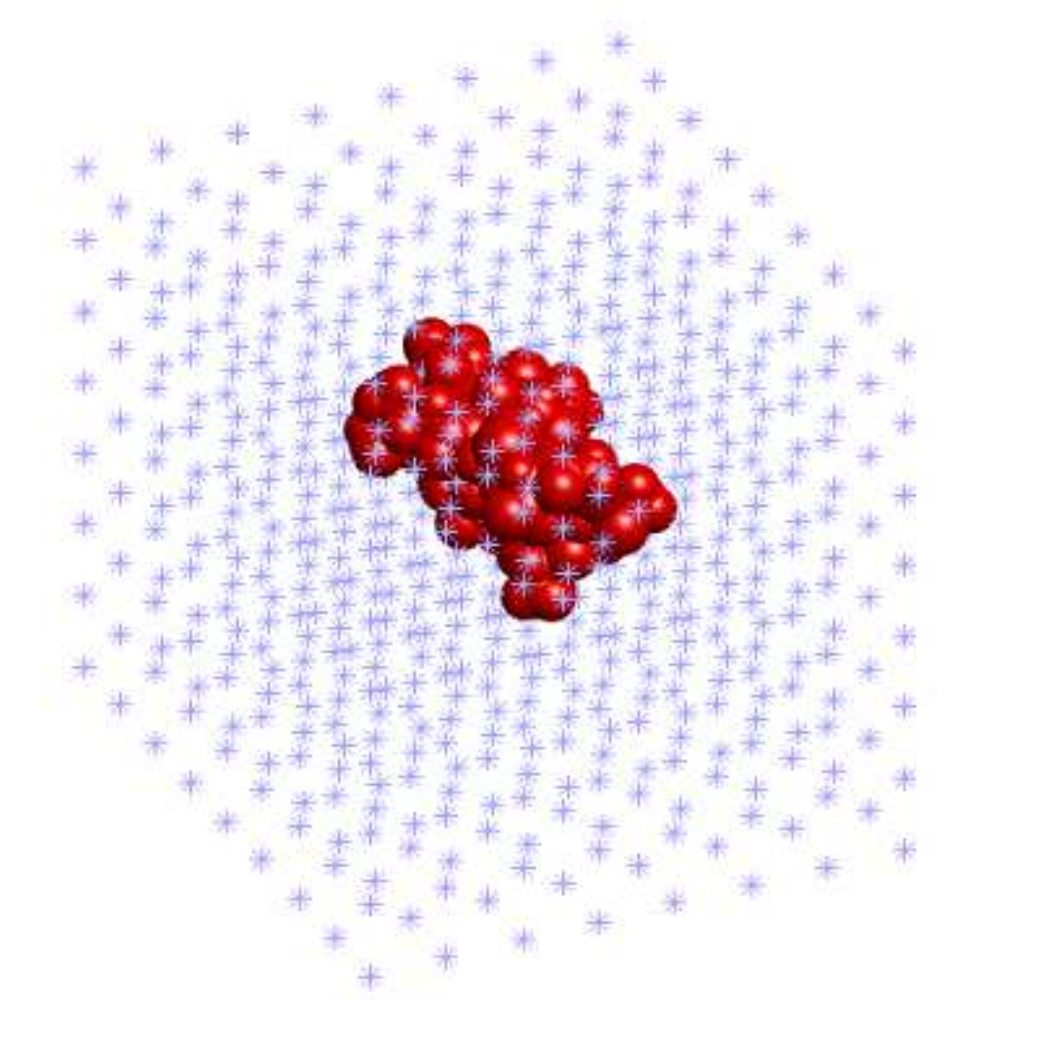}}
  \caption{training points set initialization. (a) shows a real molecule (PDBID: 1GNA) within a bounding box. (b) shows a set of uniform grid of the bounding box. (c) shows initial training points.}
  \label{point}
\end{figure}

In order to train network, in the first part, the training points set is initialized. The molecule is put in the bounding box $\Omega$ (Fig. \ref{point:surface}) in $\mathbb{R}^3$. The range of bounding box is $[a,b] \times [c,d] \times [e,f]$, where $a,b,c,d,e,f \in \mathbb{R}$. The bounding box $\Omega$ is discretized into a set of uniform grid as shown in Fig. \ref{point:box}. The training points are the grid points defined as follows
\begin{equation}
  \left\{\mathcal{P}_{ijk}\right\} = \left\{(x_{i}, y_{j}, z_{k})| x_i = a + i \cdot (b - a)/N_x, y_j = c + j \cdot (d - c)/N_y, z_k = e + k \cdot (f - e)/N_z\right\}, \\
\end{equation}
where $i=1, 2,  \cdots, N_{x}, j=1, 2, \cdots, N_{y}, k=1, 2, \cdots, N_{z}$. And $N_x, N_y, N_z$ are respectively total number of index $i, j, k$.

In the second part, the points $\left\{\mathcal{P}_{ijk}\right\}$ of train set can be labelled for training network parameter, the label of $\left\{\mathcal{P}_{ijk}\right\}$  is calculated in the following form: $\left\{\mathcal{P}_{ijk}\right\}_{label} = \phi(\left\{\mathcal{P}_{ijk}\right\})$. A set of training points $\left\{\mathbf{x}_{m}\right\}_{m=1}^{M}$ is chosen from the set of  uniform grid points $\left\{\mathcal{P}_{ijk}\right\}$, and to achieve good preservation of the molecular shape the selected points $\left\{\mathbf{x}_{m}\right\}_{m=1}^{M}$ are close to the Gaussian surface defined in Eq \ref{GaussLevelSet}. In this paper, the training points set $\left\{\mathbf{x}_{m}\right\}_{m=1}^{M}$ satisfying $\|\phi(\mathbf{x}_m)-c\|_2 \leq 1$ are selected.

\subsubsection{Parameter initialization of ellipsoid RBF neural network}
\label{Vinitialization}
In this section, we initialize the ellipsoid RBF neural network parameters $\bm{\sigma}$ defined in Eq. \ref{X}. Based on the Gaussian RBF is a degradation case of the ellipsoid Gaussian RBF, the activation function $\psi$ can be initialized as same as $\phi$.  Thus, the strategy of initialization is as follows,
\begin{enumerate}
  \item The lengths of ellipsoid RBFs neural network $\mathbf{d}$ is set to be constant vector. In this paper, the $\mathbf{d}$ can be set as follows,
  \begin{equation}\label{d_q}
    \mathbf{d} = [0.5,0.5,\cdots,0.5]^{\top} \in \mathbb{R}^{3N},
  \end{equation}
  where $N$ is the number of atom.
  \item The angles of ellipsoid RBFs neural network $\bm{\Theta}$ are set to be zeros,
  \begin{equation}\label{angle}
    \bm{\Theta} = \mathbf{0} \in \mathbb{R}^{3N}.
  \end{equation}
  \item The center coordinates $\mathbf{c}$  of ellipsoid RBFs neural network are given by the centers of atoms as follows,
  \begin{equation}\label{centercoord}
        \mathbf{c}_{i} = \left[x_{atom}^{(i)}, y_{atom}^{(i)}, z_{atom}^{(i)}\right]^{\top}, i = 1,2, \cdots,N,
  \end{equation}
  where $x_{atom}^{(i)},y_{atom}^{(i)},z_{atom}^{(i)}$ are coordinates of the $i$th atom.
  \item While $\mathbf{d}, \mathbf{c}_i, \bm{\Theta}_i$ have been chosen and atom radii $r = [r_1, r_2,\cdots ,r_N]$ is given, to initialize ellipsoid RBF activation function $\psi$ as the same as RBF $\phi$, we set the weight $\mathbf{w}$ of ellipsoid Gaussian RBF neural network as follows,
      \begin{equation}\label{weight}
        \mathbf{w} = \left[e^{\frac{r_1^2}{4}},e^{\frac{r_2^2}{4}},\cdots,e^{\frac{r_N^2}{4}}\right]^{\top}.%
      \end{equation}
\end{enumerate}

\subsubsection{Sparse optimization}
After initialization of ellipsoid RBF neural network, the sparsity of Gaussian RBF representation is computed by minimizing the loss function (Eq. \ref{Model2}). Algorithm \ref{frame} represents the main modules of our sparse optimization method, which is described below.

\begin{algorithm}
\setstretch{1.05}
\caption{Sparse optimization}
\label{frame}
\begin{algorithmic}[1]
\STATE{\textbf{Input}: PQR file including coordinates of centers and radii of atoms.}
\STATE{\textbf{Output}: The list of parameters of ellipsoid RBF neural network, i.e. solution of $\bm{\sigma}$ in minimizing
 $\mathcal{L}(\bm{\sigma})$. }
\STATE{\textbf{Step 1}. initialize network parameters $\bm{\sigma}$ as shown in Section \ref{Vinitialization}.  }
\STATE{\textbf{Step 2}. select training points set $\left\{\mathbf{x}_{m}\right\}_{m=1}^{M}$ as shown in Section \ref{Yinitialization}.}
\STATE{\textbf{Step 3}. set the number of maximum iteration $\tt MaxNiter$ and number of sparse optimization iteration $\tt SparseNiter$. set the coefficients $\rho_1$, $\rho_2$ in Eq. \ref{MODELF}.}
\STATE{\textbf{Step 4}. initialize the variable of iteration: $\tt Niter= 0$ and set tolerance: $\tt tol_1 $ and $\tt tol_2$. }
\STATE{\textbf{Step 5}. select the size of batch $\tt Batch\_size$ for optimization algorithm.}
\STATE{\textbf{Step 6}. optimize two terms of loss function in Eq. \ref{Model2} alternatively.}
\WHILE{$\tt Niter < \tt MaxNiter$}
\STATE {$\tt Niter = \tt Niter+ 1$}
\STATE \textbf{Step 6.1}. prune the useless the ellipsoid RBF neuron $|w_i| < \tt tol_1$ every $\tt check_{step}$ steps.
\STATE \textbf{Step 6.2}. calculate $\psi(\mathbf{x}_{m})$ for training points set by $\bm{\sigma}$.
\STATE \textbf{Step 6.3}. check the maximum of error between $\psi$ and $\phi$ at training points set $\left\{\mathbf{x}_m\right\}_{m=1}^M$ and update the coefficients $\rho _s$ and $\rho _{l}$.
\IF {$\mathop {\max } _{1 \leq m\leq M} \| \psi(\mathbf{x}_m) - \phi(\mathbf{x}_m) \| > \tt tol_2$}
\STATE $\quad\quad$ $\rho_1 = 1, \rho_2 = 0$
\ENDIF
\STATE \textbf{Step 6.4}. accuracy optimization for $E_s$ by set coefficients $\rho_1$ and $\rho_{l}$.
\IF {$\tt Niter > SparseNiter$}
\STATE $\quad\quad$ $\rho_1 = 1, \rho_2 = 0$
\ENDIF
\STATE \textbf{Step 6.5}. optimize the loss function $\mathcal{L}(\bm{\sigma})$ by batch ADAM algorithm.
\ENDWHILE
\end{algorithmic}
\end{algorithm}

Step 1 shows initialization of parameters $\bm{\sigma}$ for ellipsoid RBF neural network. Step 2 selects training points set $\left\{\mathbf{x}_{m}\right\}_{m=1}^{M}$. Step 3 and step 4 initialize some variables, i.e., the number of total iterations, the number of sparse optimization iterations, error tolerance and the coefficients $\rho_1$, $\rho_2$ . Step 5 sets the size of batch ($\tt Batch\_size = 1000$) for optimization algorithm. Step 6 shows the numerical algorithm of optimization for loss function (Eq. \ref{Model2}). Step 6.1 prunes useless the ellipsoid RBF neuron if the corresponding weight $w_i$ connecting the $i$th hidden neuron to the output layer is less than $\tt tol_1$ per $\tt check_{step}$ steps. In this paper, we set $\tt tol_1 = \tt 1e-3$ and $\tt check_{step} = 20$. Step 6.2 calculates the prediction value $\psi$ for all training points set. Step 6.3 checks the maximum of error between $\psi$ and $\phi$ at training points set $\left\{\mathbf{x}_{m}\right\}_{m=1}^{M}$ and updates the coefficients $\rho _s$ and $\rho _{l}$, where $\tt tol_2 = 0.1$. Step 6.4, after doing $\tt SparseNiter$ iterations, with the number of effective neurons fixed, keeps doing some steps of minimization of $E_s$ to achieve better accuracy of the approximation on training points set. Step 6.5 updates the network parameter $\bm{\sigma}$ and optimizes loss function $\mathcal{L}(\bm{\sigma})$ by batch ADAM method. The pipeline of step 6.5 is as follows,
\begin{equation}\label{ADMM}
  \bm{\sigma}_{k+1} =  \bm{\sigma}_{k} - \tau \frac{\frac{\beta_1}{1 - \beta_1^k}\cdot m_{k-1} + \frac{1 - \beta_1}{1 - \beta_1^k}\cdot \nabla \mathcal{L}_k}{\sqrt{\frac{\beta_2}{1 - \beta_2^k}\cdot v_{k-1}+\frac{1 - \beta_2}{1 - \beta_2^k}\cdot (\nabla \mathcal{L}_k)^2} + \epsilon},
\end{equation}
where $\tau=0.002$ is learning rate. $\beta_1$, $\beta_2$ and $\epsilon$ are set for default value ($\beta_1 = 0.9, \beta_2 = 0.999, \epsilon = 10^{-8}$), $m_{k}$ is $k$th biased first moment estimate. $v_{k}$ is the $k$th biased second raw estimate.



\section{Results and discussion}\label{section3}
In this section, we present some numerical experimental examples to illustrate the effectiveness of our network and method for representing the Gaussian surface sparsely. Comparisons are made among our network, the original definition of  Gaussian density maps and sparse RBF method \cite{2016LiSparse}. A set of biomolecules taken from the RCSB Protein Data Bank is chosen as a benchmark set. The number of atoms in these biomolecules ranges from hundreds to thousands. These molecules are chosen randomly from RCSB Protein Data Bank, and no particular structure is specified.  The implementation of the algorithms is based on the PyTorch. All computations were run on an Nvidia Tesla P40 GPU. Further quantitative analysis of the result is given in the following subsections.

\subsection{Sparse optimization results}
Twenty biomolecules are chosen to be sparsely represented by the ellipsoid RBFs neural network and sparse RBF method \cite{2016LiSparse}.  For fair comparison, the initial centers of RBFs are selected to be atom center coordinates for both methods. Table \ref{ProteinResult} shows the final number of effective basis from the results of our method and sparse RBF method.
\begin{table}[H]
\begin{center}
\caption{The number of atoms for 20 test proteins. The fourth  column shows the number of RBFs by sparse RBF method. The last  column shows the number of the ellipsoid RBF neural network. The decay rate $d$ in Eq. \ref{GaussKernel} is uniformly taken as 0.5.}\label{ProteinResult}
\begin{tabular}{cc|ccc}
\toprule[2pt]
INDEX & PDBID & NATOM & Sparse RBF & \textbf{Our method} \\  \midrule[1pt]
1 & ADP & 39 & 13 & \textbf{5} \\ 
2 & 2LWC & 75 & 51 & \textbf{6} \\ 
3 & 3SGS & 94 & 56 & \textbf{9} \\ 
4 & 1GNA & 163 & 108 & \textbf{17} \\ 
5 & 1V4Z & 266 & 266 & \textbf{22} \\ 
6 & 1BTQ & 307 & 252 & \textbf{25} \\ 
7 & 2FLY & 355 & 267 & \textbf{28} \\ 
8 & 6BST & 478 & 315 & \textbf{49} \\ 
9 & 1MAG & 552 & 502 & \textbf{46} \\ 
10 & 2JM0 & 569 & 424 &  \textbf{52} \\ 
11 & 1BWX & 643 & 537 & \textbf{54} \\ 
12 & 2O3M & 714 & 566 & \textbf{62} \\ 
13 & FAS2 & 906 & 722 & \textbf{76} \\ 
14 & 2IJI & 929 & 742 & \textbf{72} \\ 
15 & 3SJ4 & 1283 & 953 & \textbf{132} \\ 
16 & 3LOD & 2315 & 1810 & \textbf{180} \\ 
17 & 1RMP & 3514 & 2871 & \textbf{271} \\ 
18 & 1IF4 & 4251 & 3288 & \textbf{301} \\ 
19 & 1BL8 & 5892 & 3491 & \textbf{452} \\ 
20 & AChE & 8280 & 4438 & \textbf{748} \\ 
\bottomrule[2pt]
\end{tabular}
\end{center}
\end{table}

\begin{figure}[H]
  \centering
  \includegraphics[scale = 0.65]{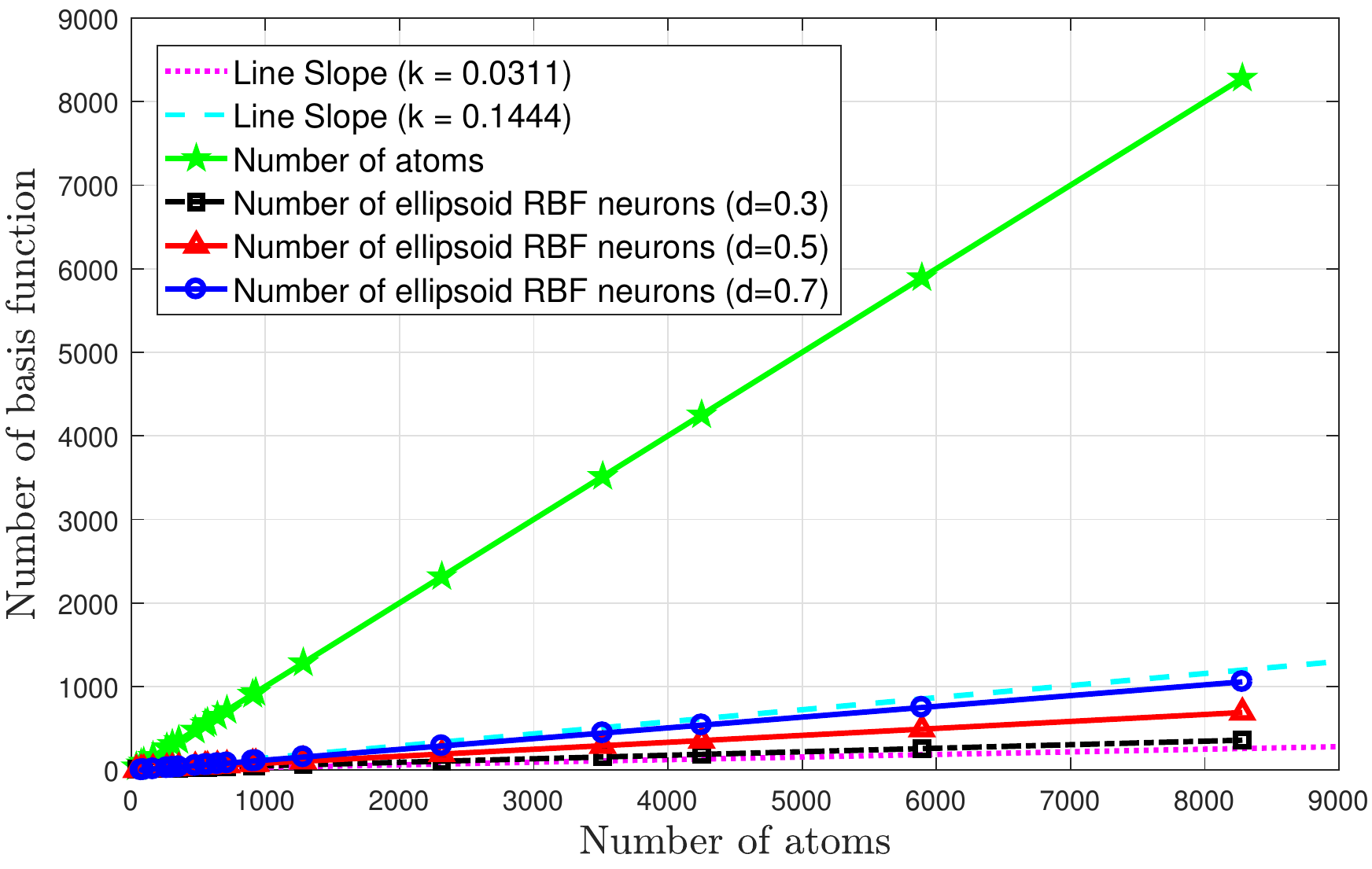}\\
  \caption{Relationship between the number of atoms and the number of ellipsoid RBF neurons after sparse representation.}\label{compressResult}
\end{figure}

Fig. \ref{compressResult} presents the relation between the number of ellipsoid RBF neurons in final sparse representation and the number of atoms in the corresponding molecule. The number of atoms for the original Gaussian molecular surfaces is shown by green lines with pentagram markers. To present sparsity of final results from our method, we define the sparse ratio $R_s$ as: $R_s = \frac{N_{ERBF}}{N_{ATOM}}$, where $N_{ERBF}$ is the number of ellipsoid RBF neurons and $N_{ATOM}$ presents the number of atoms. In Fig. \ref{compressResult}, the changes of sparse ratios with respect to number of atoms for different decay rates ($d$ in Eq. \ref{GaussKernel} equals to 0.3, 0.5 and 0.7) are shown by solid lines with square, circle and triangle markers. The slope of dashed line is the lower bound of sparse ratio ($k=  0.0311$). The slope of dotted line is the upper bound of sparse ratio ($k= 0.1444$). The sparse ratios in the results of our numerical experiments are in ($0.0311, 0.1444$). The results show that the larger the decay rate $d$ (leading to more rugged and complex molecular surface), the bigger the sparse ratio is going to be. The sparse ratios for the Gaussian molecular surface with $d = 0.3$ are smaller than those of Gaussian molecular surface with $d = 0.7$ as shown in Fig. \ref{compressResult}.

\begin{figure}[H]
  \centering
    \subfigure[]{
  \label{CostNERBF:lossfunction}
  \includegraphics[scale = 0.5]{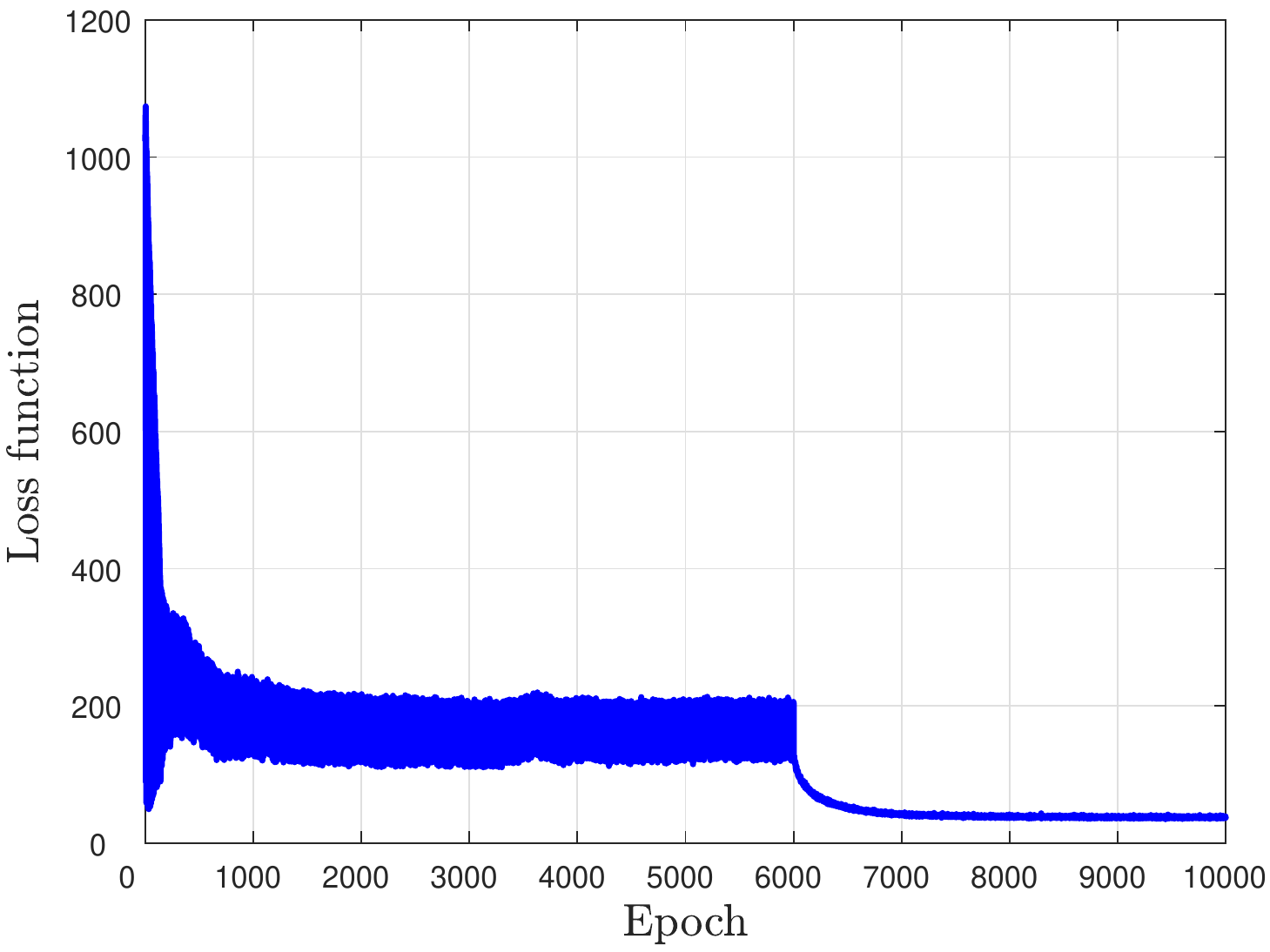}}
  \subfigure[]{
  \label{CostNERBF:sizeRBF}
  \includegraphics[scale = 0.5]{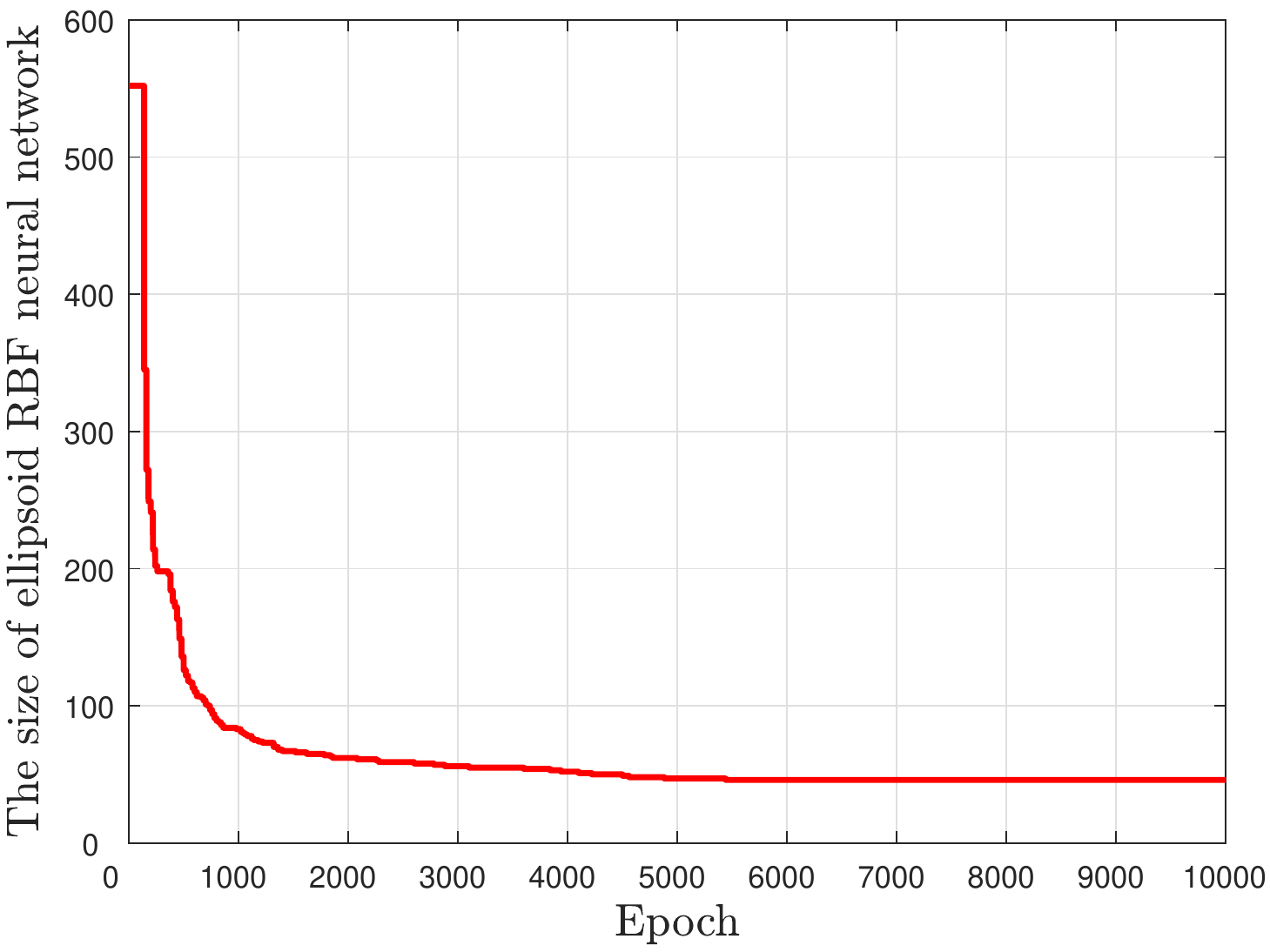}}
  \caption{One test of our algorithm on molecule 1MAG. The blue curve is the objective function trajectory during the 10000 iterations. The red line represents the number of basis functions. The number of initial ellipsoid RBF neurons for this trial is 552 and the number of final ellipsoid RBF neurons is 46.}\label{CostNERBF}
\end{figure}

Fig. \ref{CostNERBF} shows the loss function and the number of ellipsoid RBF neurons is decreasing as the number of iterations increases in the experiment for molecule 1MAG. In this experiment, the $\tt MaxNiter$ and $\tt SparseNiter$ are set to be 10000 and 6000, respectively. After 6000 iterations, $\rho_1$ is set to be zero to minimize $E_s$ term  solely, thus the value of loss function has an abrupt change. The number of ellipsoid RBF neurons decreases dramatically during the iteration process. As shown in Fig.  \ref{CostNERBF}, the model with 46 ellipsoid RBF neurons achieves the minimum error with a relatively less number of ellipsoid RBF neurons.

\begin{figure}[H]
  \centering
  \includegraphics[scale = 0.6]{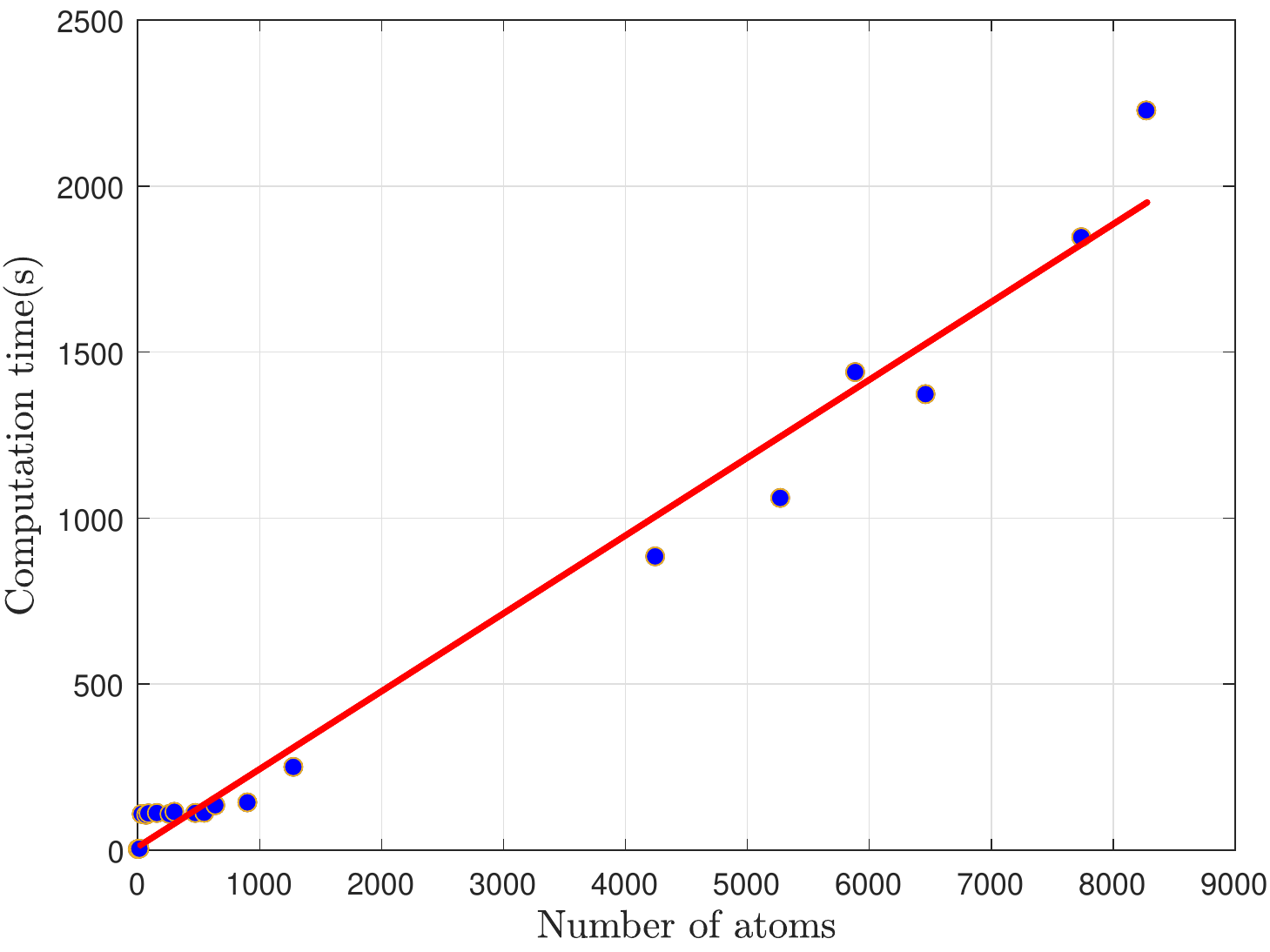}\\
  \caption{Computational performance of training algorithm for our network.}\label{CPUTIME}
\end{figure}

The complexity of training algorithm for our network is almost $O(N)$, which is shown in Fig \ref{CPUTIME}.

\subsection{Shape preservation and further results analysis}
In this subsection, we first check whether the Gaussian surface is preserved after the process of sparse representation through our method. The area, the enclosed volume and the Hausdorff distance are the three criteria to judge whether two surfaces are close enough. These criteria can be calculated on the triangular mesh of the surface. The triangular meshes of molecular surfaces before and after sparse representation  are generated through $\tt isosurface$ function in MATLAB. For a triangular surface mesh, the surface area $S$ is determined using the following equation:
\begin{equation}
S=\frac{1}{2} \sum_{i=1}^{n_f}\left\|\overrightarrow {V_{1}^{i} V_{2}^{i}} \times \overrightarrow{V_{1}^{i} V_{3}^{i}}\right\|,
\end{equation}
where $n_f$ is the number of triangle elements and $V_{1}^{i}, V_{2}^{i}, V_{3}^{i}$ denote the coordinates of the three vertices for the $i$th triangle.

The volume V enclosed by the surface mesh is determined using the following equation:
\begin{equation}
V=\frac{1}{6} \sum_{i=1}^{n_f} \overrightarrow{V_{2}^{i} V_{1}^{i}} \times \overrightarrow{V_{3}^{i} V_{1}^{i}} \bullet \vec{c}_{i},
\end{equation}
where $c_i$ is the vector from the center of the $i$th triangle to the origin.

The relative errors of area/volume and the Hausdorff distance are used to characterize the difference between the surfaces before and after sparse representation. The relative errors of area and volume are calculated using the following formulas:
\begin{equation}
Error_A = \frac{|A_{our} - A_{original}|}{A_{original}},
\end{equation}
\begin{equation}
Error_V = \frac{|V_{our} - V_{original}|}{V_{original}},
\end{equation}
where $A_{original}$ and $A_{our}$ denote the surface areas of the original and our sparsely represented surfaces respectively. $V_{original}$ and $V_{our}$ denote the corresponding surface volumes of the original and our surfaces respectively.

The Hausdorff distance between two surface meshes is defined as follows,
\begin{equation}
H(S_1,S_2) = max\left(\mathop{\textup{max}}_{p \in S_1}e(p,S_2),\mathop{\textup{max}}_{p \in S_2}e(p,S_1)\right),
\end{equation}
where
\begin{equation}
e(p,S) = \mathop{\textup{min}}_{p' \in S}d(p,p'),
\end{equation}
$S_1$ and $S_2$ are two piecewise surfaces spanned by the two corresponding meshes, and $d\left(p, p^{\prime}\right)$ is the Euclidean distance between the points $p$ and $p^{\prime}$. In our work, we use Metro \cite{Cignoni1996MME} to compute the Hausdorff distance.

The areas and the volumes enclosed by the surface before and after the sparse representation for each of the molecules are listed in Table \ref{tab:geo}. The Hausdorff distances between the original surface and the final surface for the biomolecules are also listed in Table \ref{tab:geo}.

\begin{table}[H]
\begin{center}
\caption{The areas, volumes and Hausdorff distances obtained with the original and the final surfaces for ten biomolecules. Note: isovalue $\phi = 1.0$, initial decay rate $d = 0.5$.}\label{tab:geo}
\begin{tabular}{c|c|c|c|c|c|c|c}
\hline \hline
\multirow{2}*{Molecule}  & \multicolumn{3}{c}{Area ($\AA^2$)} & \multicolumn{3}{|c|}{Volume ($\AA^3$)} & \multirow{2}*{Distance ($\AA$)} \\
\cline{2-7}  
 & Original & Our & $Error_A$ & Original & Our & $Error_V$ &   \\ \hline
ADP	&367.9334 &	358.4047& 	0.0259& 	458.0317& 	454.5578& 	0.0076& 	0.6605 \\ 
2LWC&	504.8863 &	494.5004& 	0.0206& 	856.9564 &	850.5393& 	0.0075& 	0.4497 \\ 
1GNA&	1006.1213& 	995.4826& 	0.0106& 	1862.7883& 	1855.4815 &	0.0039& 	0.4764 \\ 
1BTQ&	1782.1843& 	1749.4445 &	0.0184 	&3412.7345& 	3406.8652& 	0.0017& 	0.6027 \\ 
1MAG&	2479.4398& 	2438.4246& 	0.0165& 	5732.9858& 	5700.8069& 	0.0056 &	0.5441 \\ 
1BWX&	2925.0557& 	2890.7706& 	0.0117& 	6638.2112& 	6609.3813& 	0.0043& 	0.7311 \\ 
FAS2&	3771.6093 &	3690.4698 &	0.0215& 	9198.0803 &	9168.8722 &	0.0032& 	0.7484 \\ 
2IJI	&3783.6502 &	3731.7199 	&0.0137& 	9537.9781 &	9502.8469 &	0.0037& 	0.6187 \\ 
3SJ4	&5887.9106 &	5797.6176 &	0.0153& 	13208.3953& 	13175.7877& 	0.0025& 	0.7209 \\ \hline \hline
\end{tabular}
\end{center}
\end{table}

\begin{figure}
  \centering
  \includegraphics[scale = 0.45]{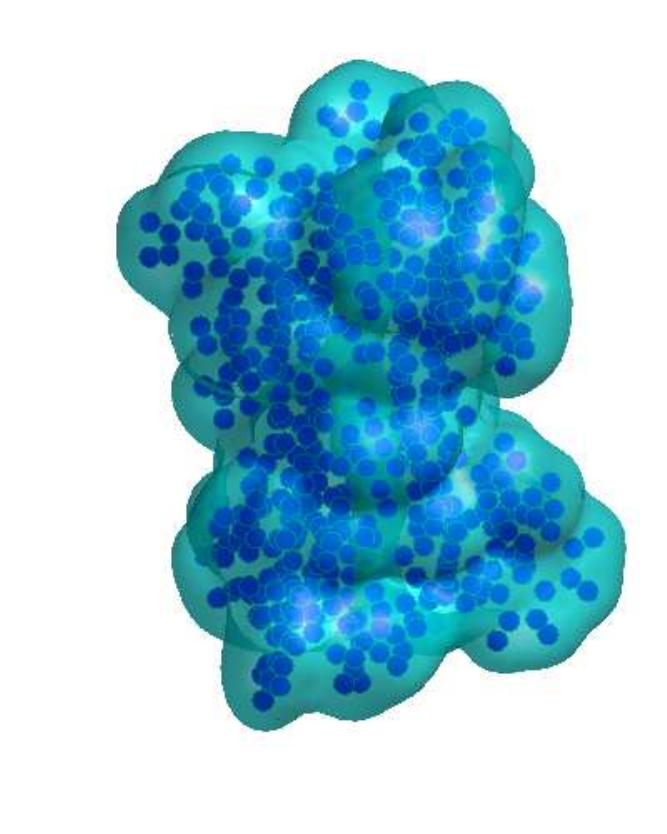}
  \includegraphics[scale = 0.45]{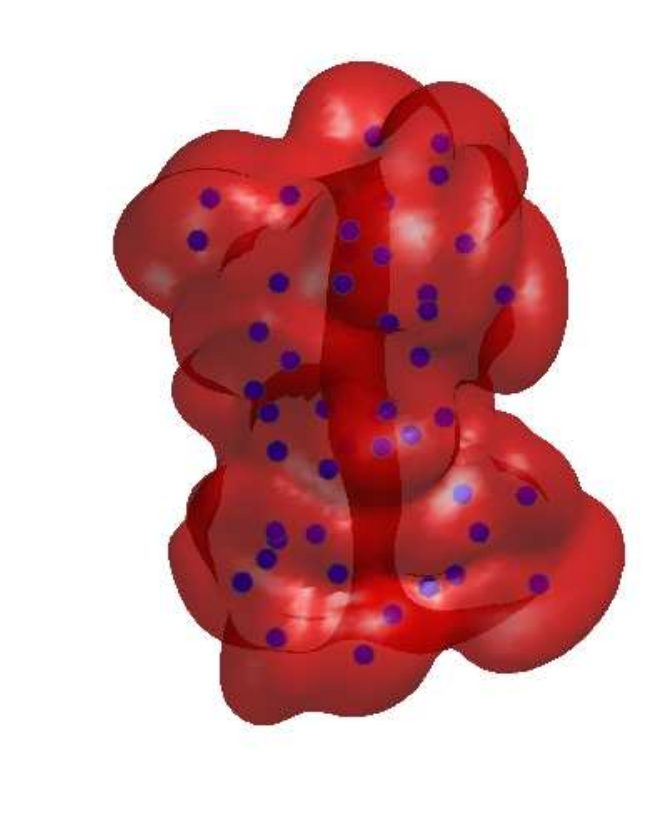}
  \includegraphics[scale = 0.45]{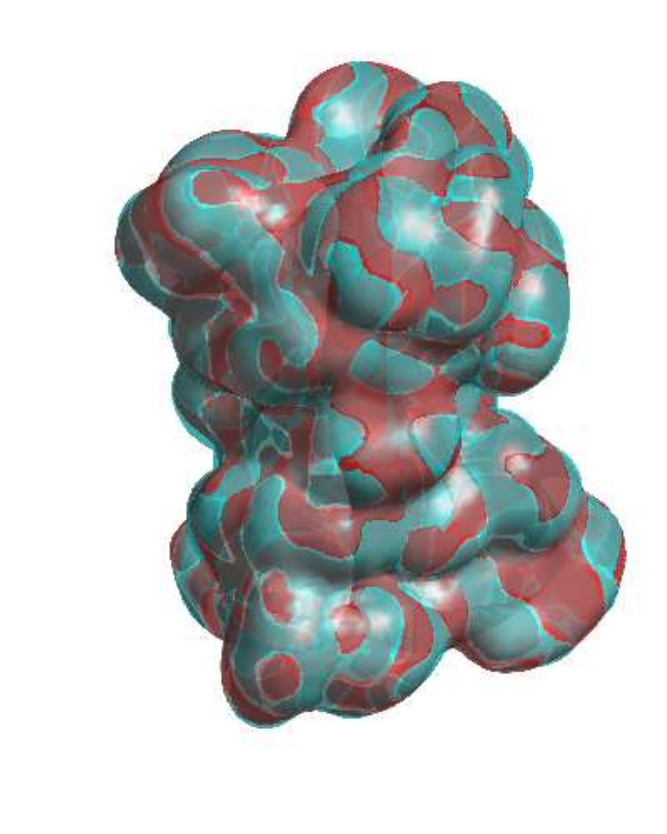}
  \includegraphics[scale = 0.45]{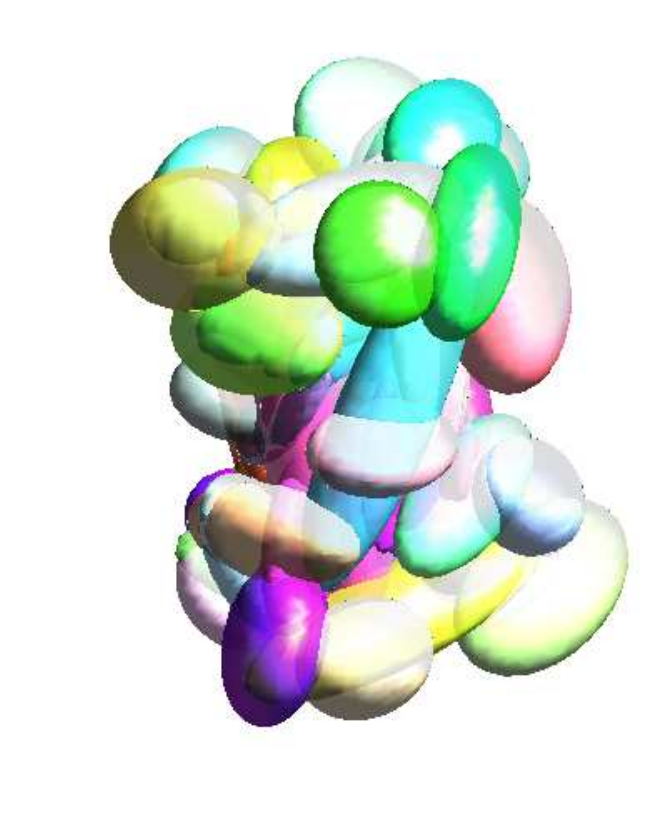}\\
  \includegraphics[scale = 0.45]{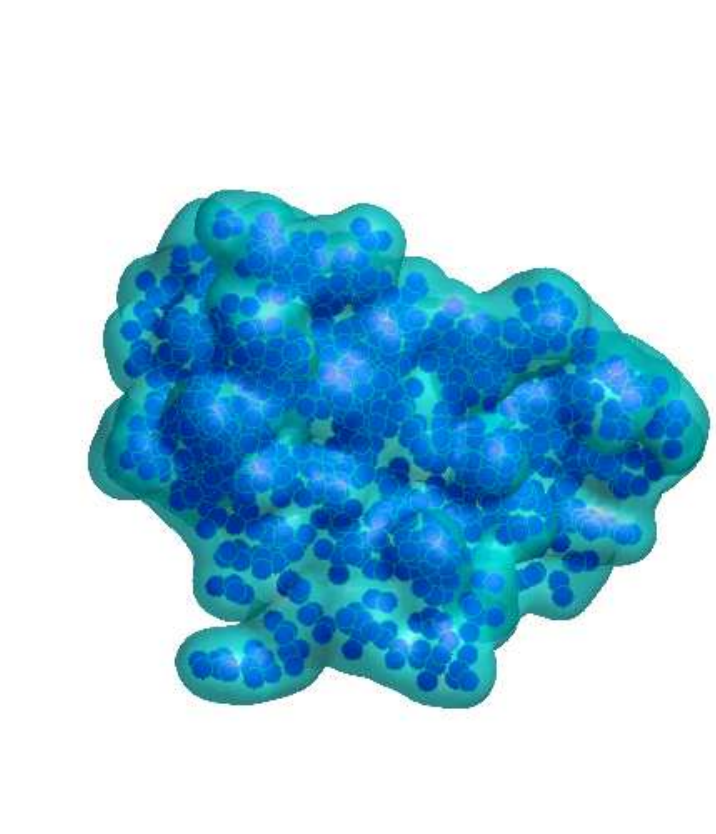}
  \includegraphics[scale = 0.45]{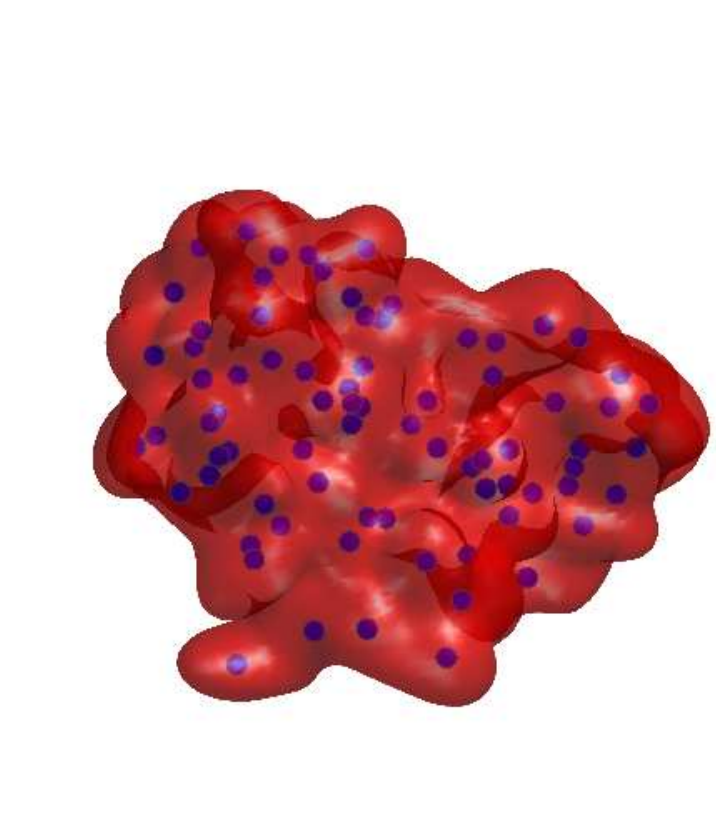}
  \includegraphics[scale = 0.45]{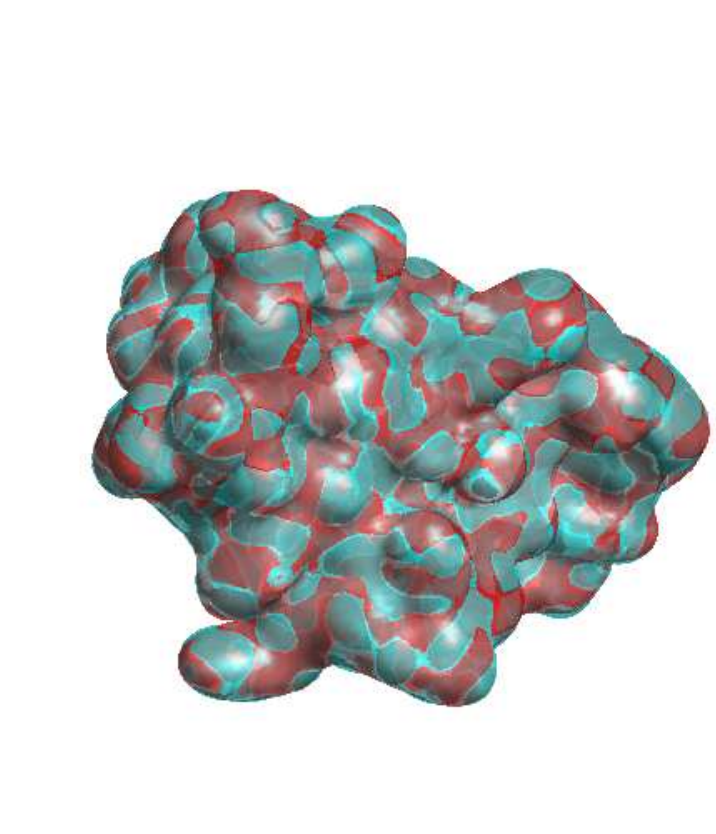}
  \includegraphics[scale = 0.45]{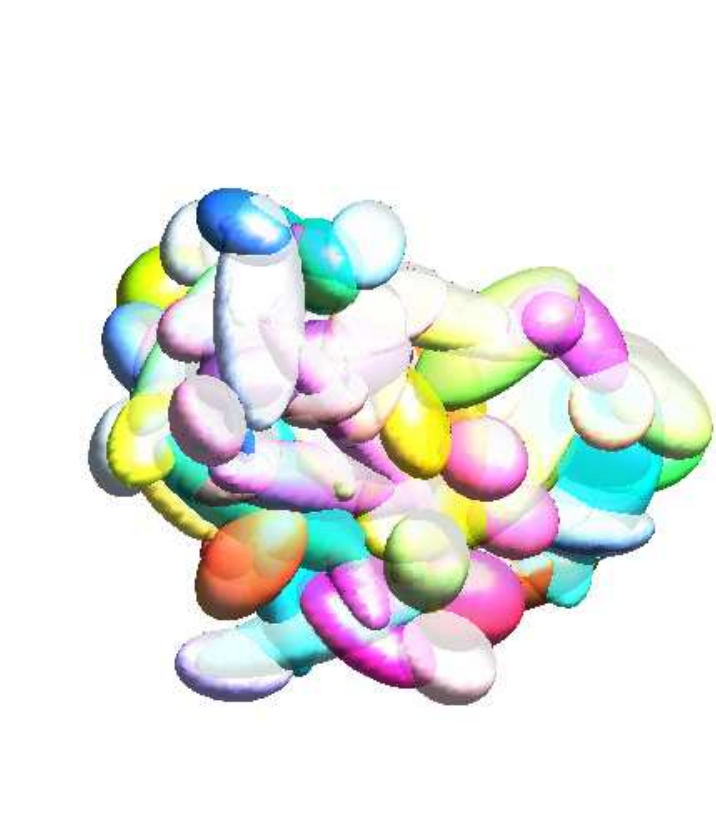}\\
  \includegraphics[scale = 0.45]{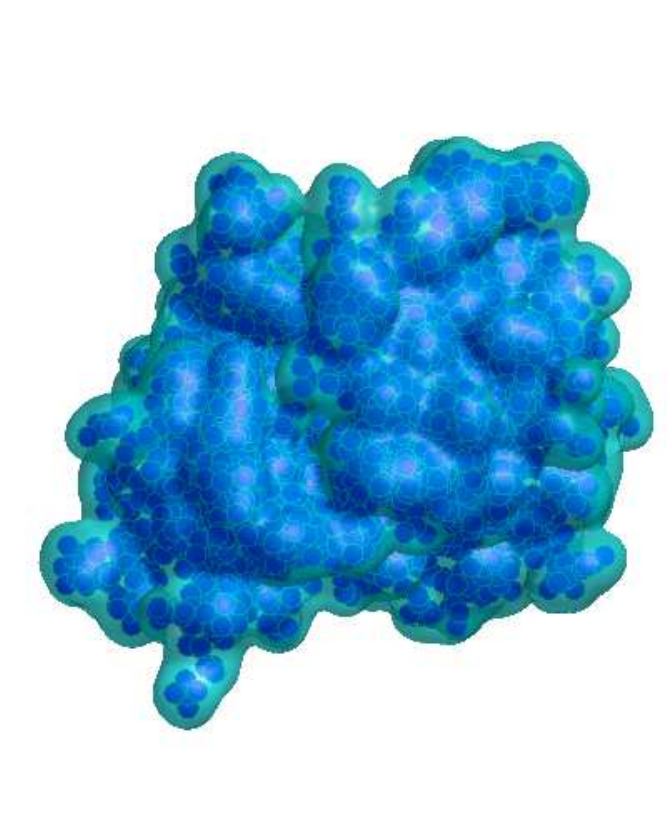}
  \includegraphics[scale = 0.45]{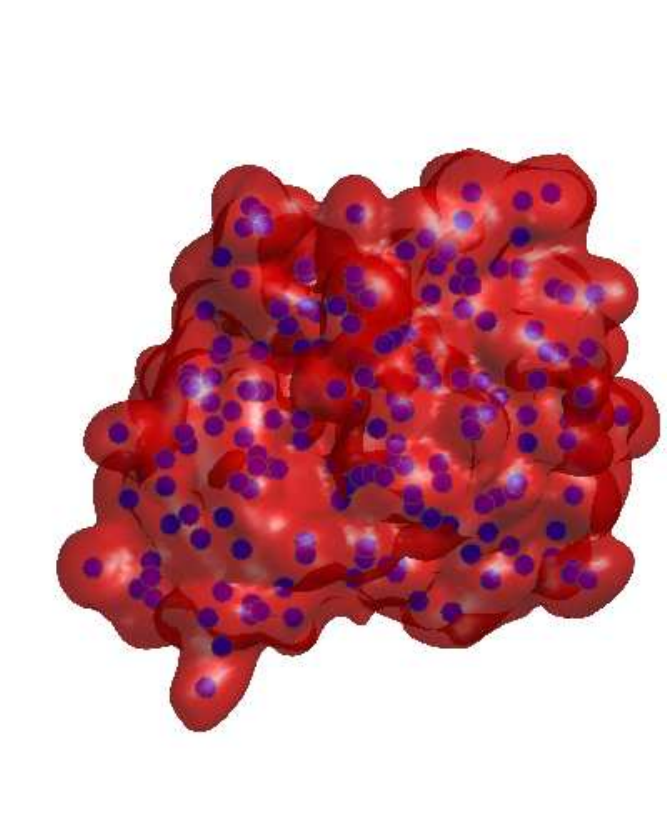}
  \includegraphics[scale = 0.45]{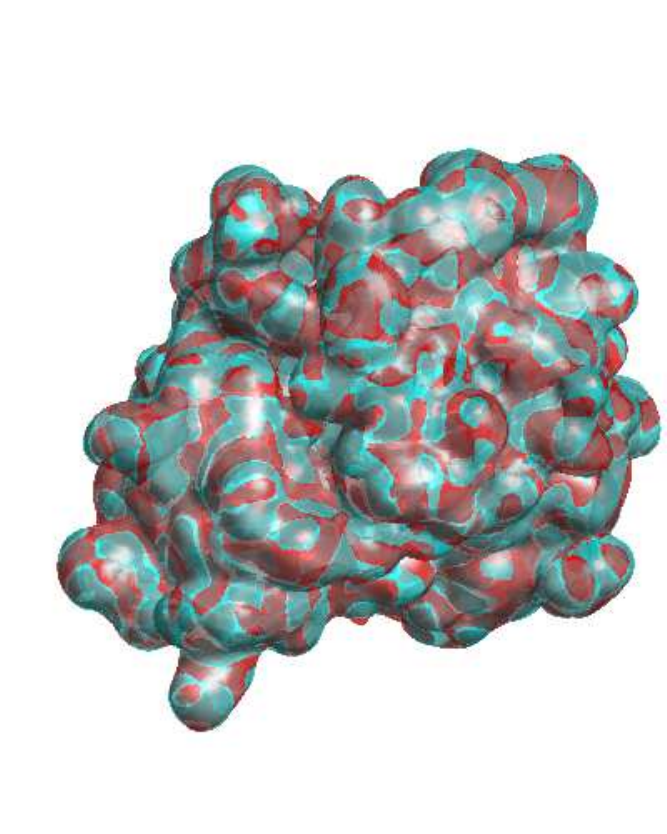}
  \includegraphics[scale = 0.45]{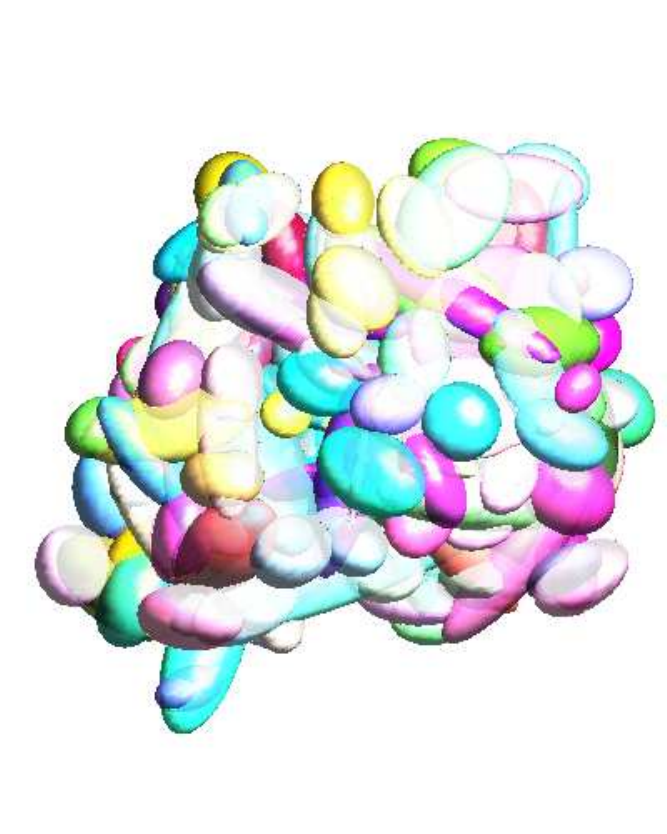}\\
  \includegraphics[scale = 0.45]{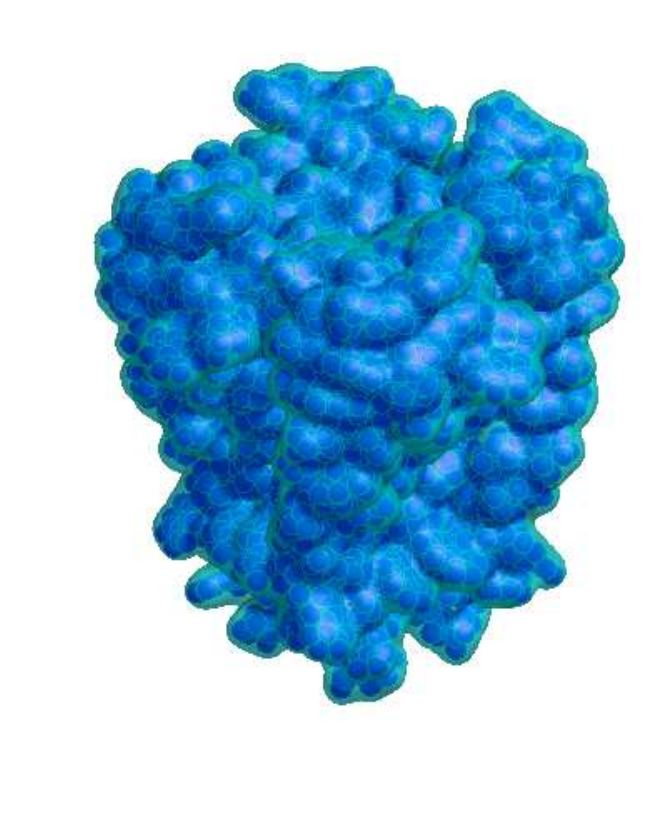}
  \includegraphics[scale = 0.45]{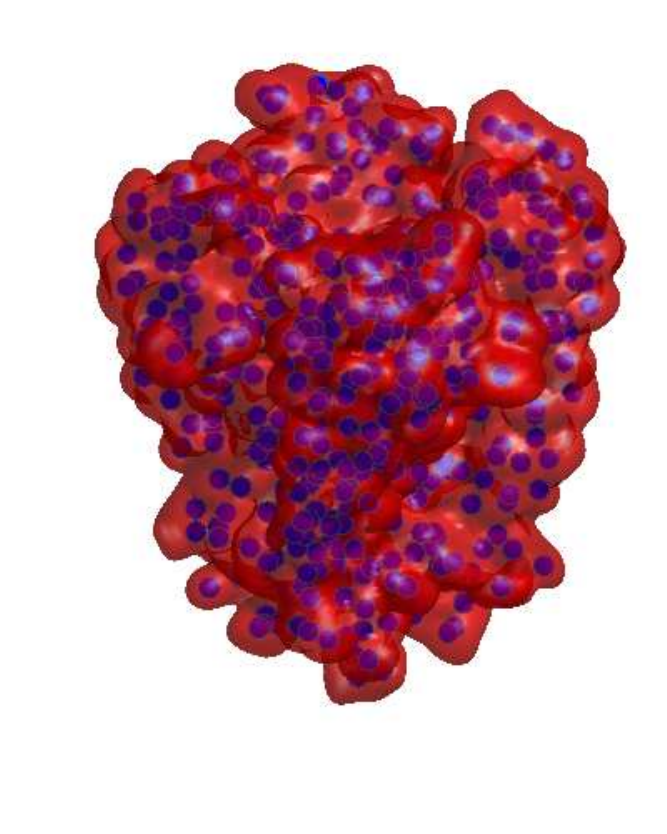}
  \includegraphics[scale = 0.45]{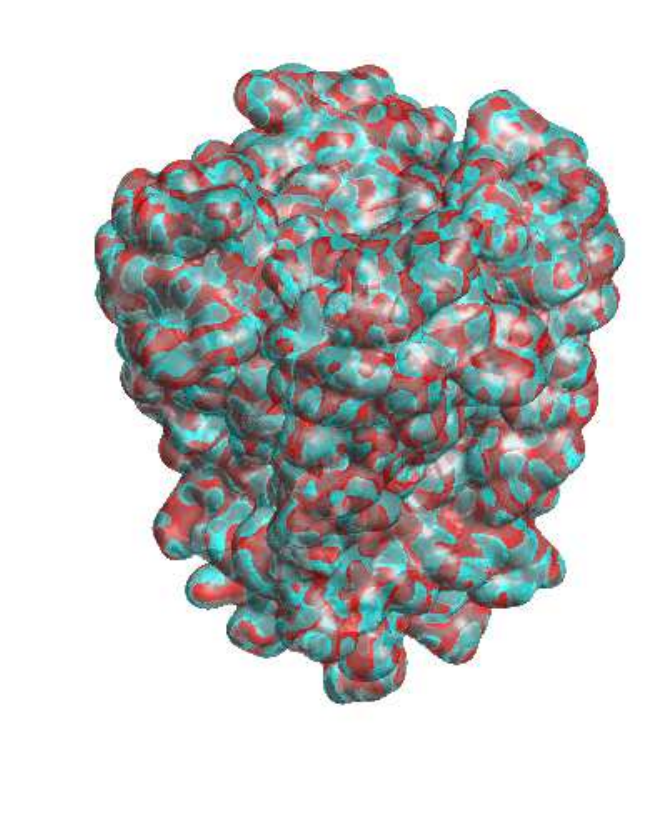}
  \includegraphics[scale = 0.45]{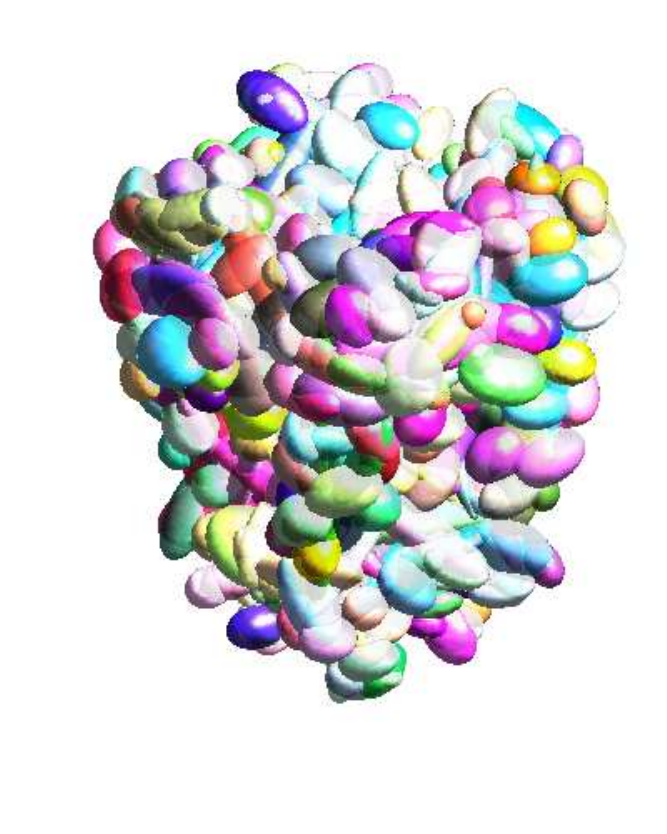}\\
  \includegraphics[scale = 0.45]{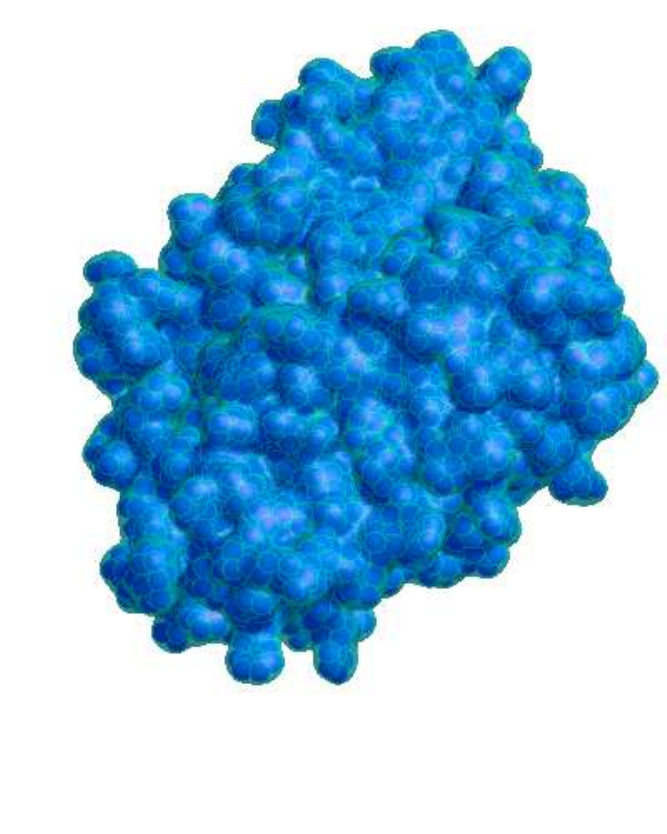}
  \includegraphics[scale = 0.45]{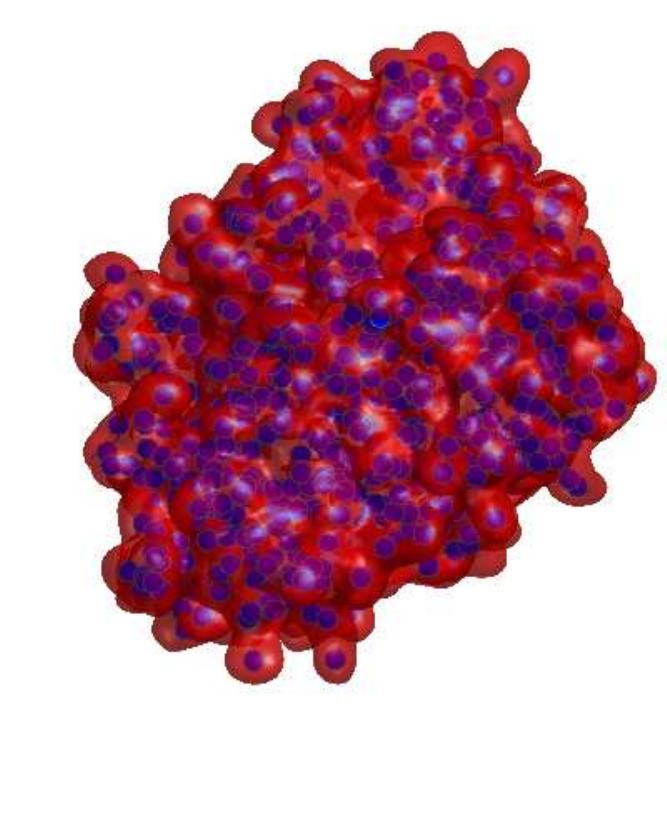}
  \includegraphics[scale = 0.45]{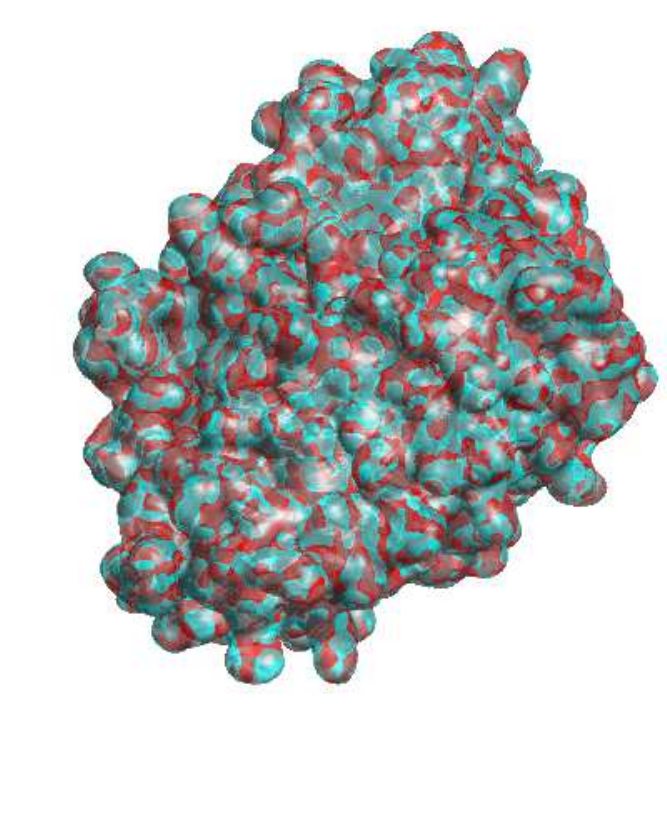}
  \includegraphics[scale = 0.45]{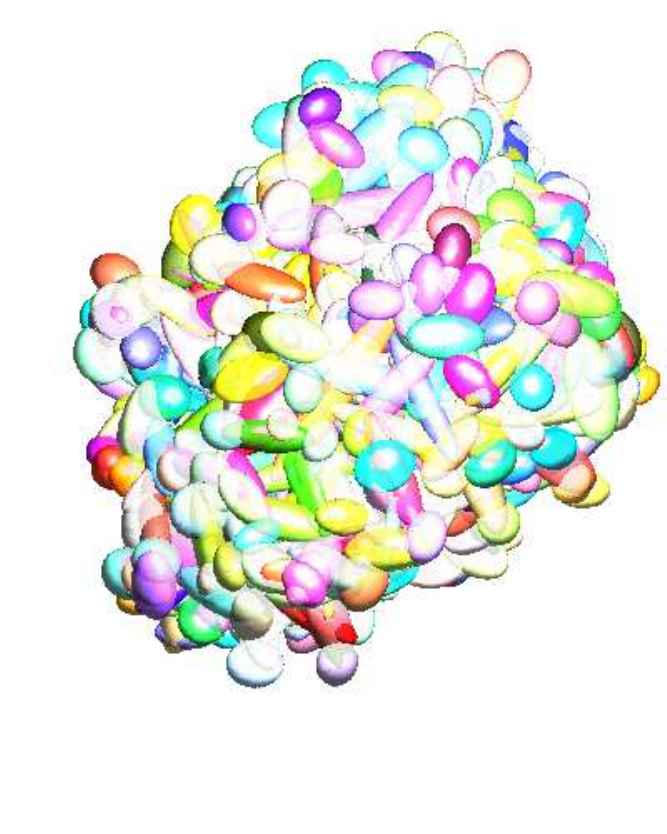}
  \caption{Fitting results of our optimization algorithm. Left to right: Original surface (left column), Final surface (middle left column) and Original surface overlapped with Final surface (middle right column),
   the ellipsoid Gaussian RBFs in the sparse representation from our method (last column). From top to bottom: 1MAG (first row), FAS2 (second row), 3LOD (third row), 1BL8 (fourth row) and AChE (fifth row). The blue points represent the locations of Gaussian RBF centers.}
  \label{mixresult}
\end{figure}

Fig. \ref{mixresult} illustrates some fitted results of the sparse optimization model. The first column shows original Gaussian surface for five molecules. The second column is the final ellipsoid Gaussian surface in our method, where the blue points represent the location of Gaussian RBF centers. It indicates our method needs less number of ellipsoid RBFs neurons to represent surface. The third column is the original surface overlapped with the final surface. It shows that the final surface is close to the original surface.
The last column shows the configurations of ellipsoid RBF neurons in the sparse representation of five molecules from our method. It demonstrates that after the process of sparse representation, the number of ellipsoid RBF neurons are much sparser than the RBFs in the original definition of Gaussian surface. And, obviously, each ellipsoid RBF is a local shape descriptor of the molecular shape.

\subsection{Electrostatic solvation energy calculation based on the sparsely represented surface}
\begin{figure}[H]
  \centering
  \includegraphics[scale = 0.18]{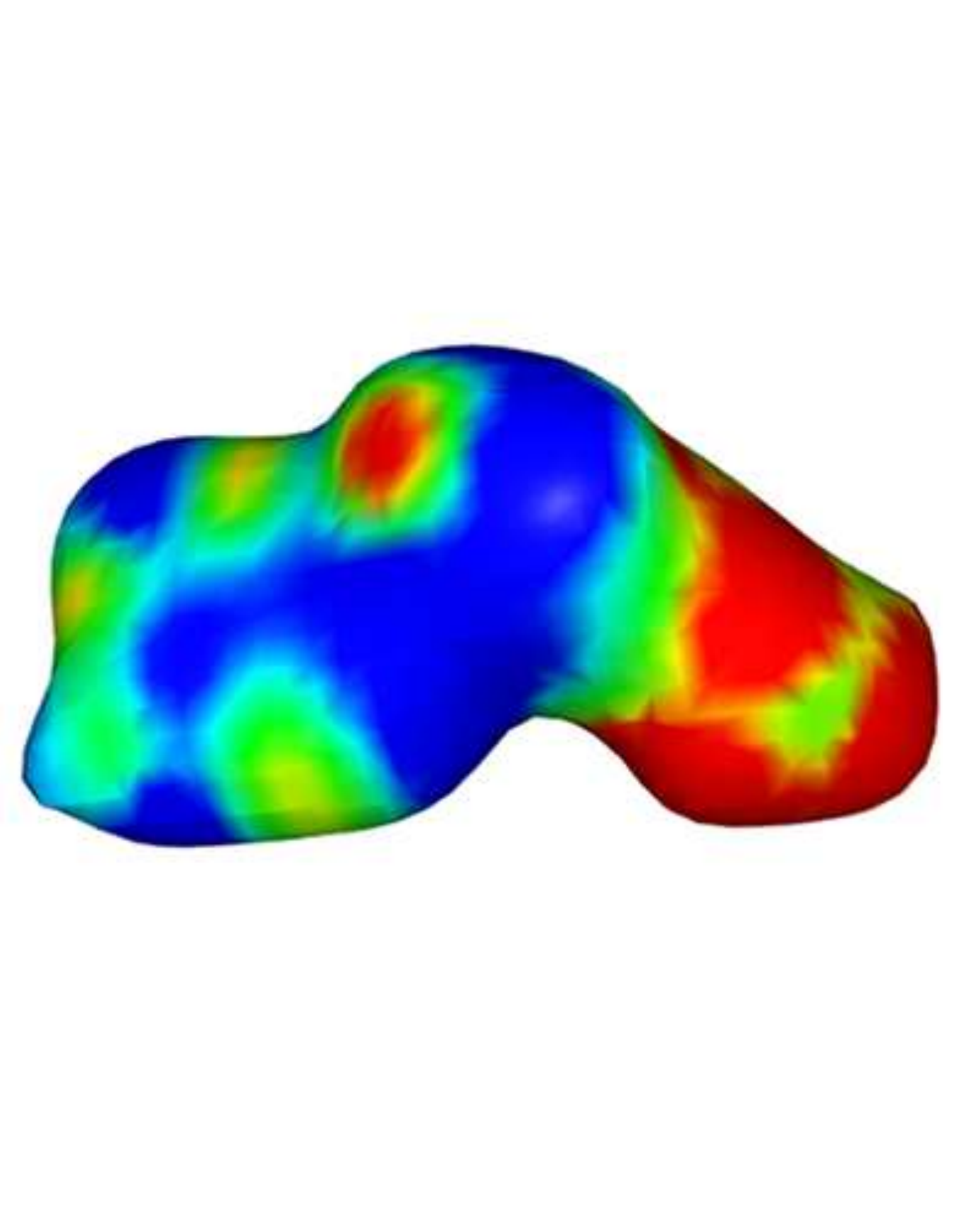}
  \includegraphics[scale = 0.18]{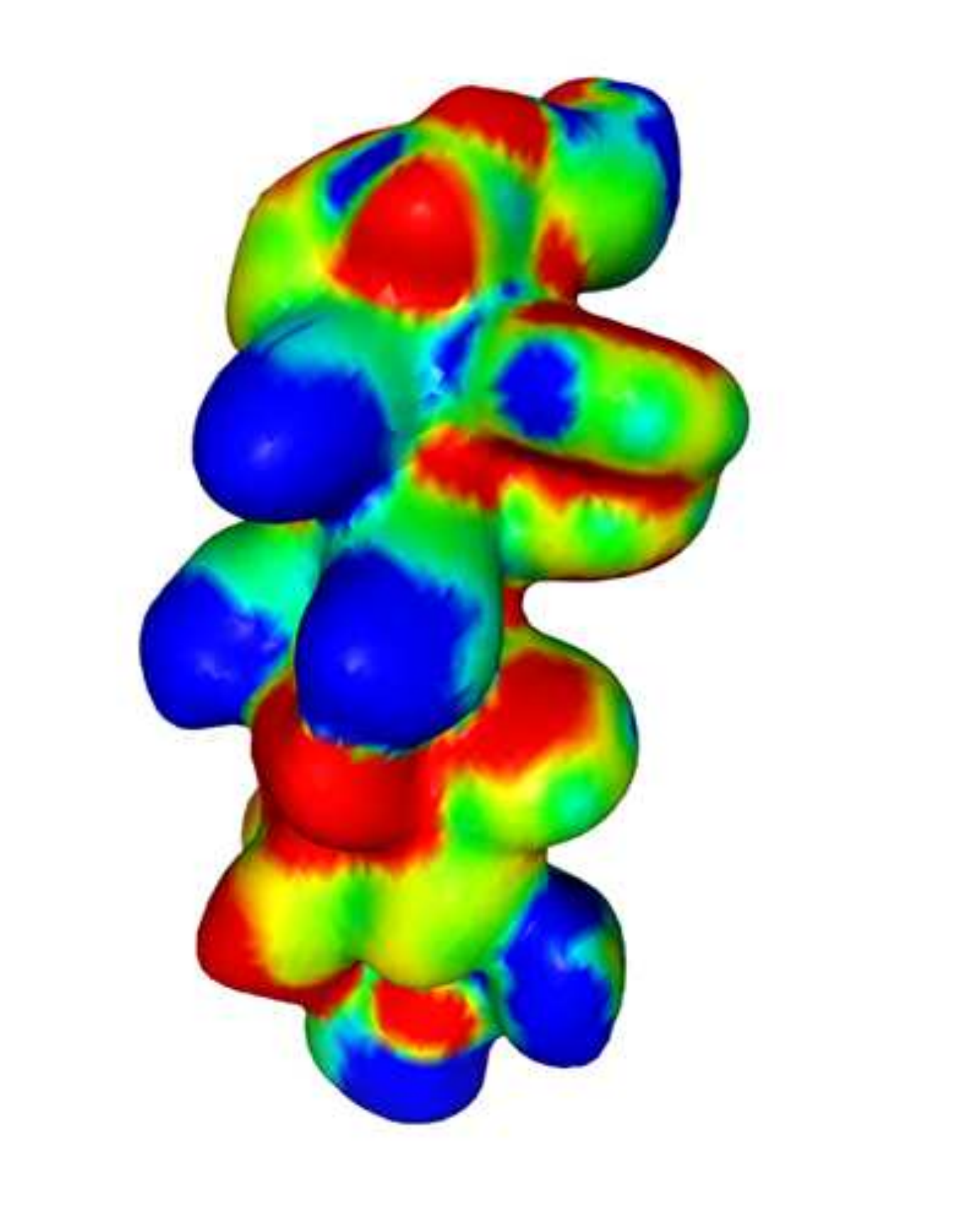}
  \includegraphics[scale = 0.18]{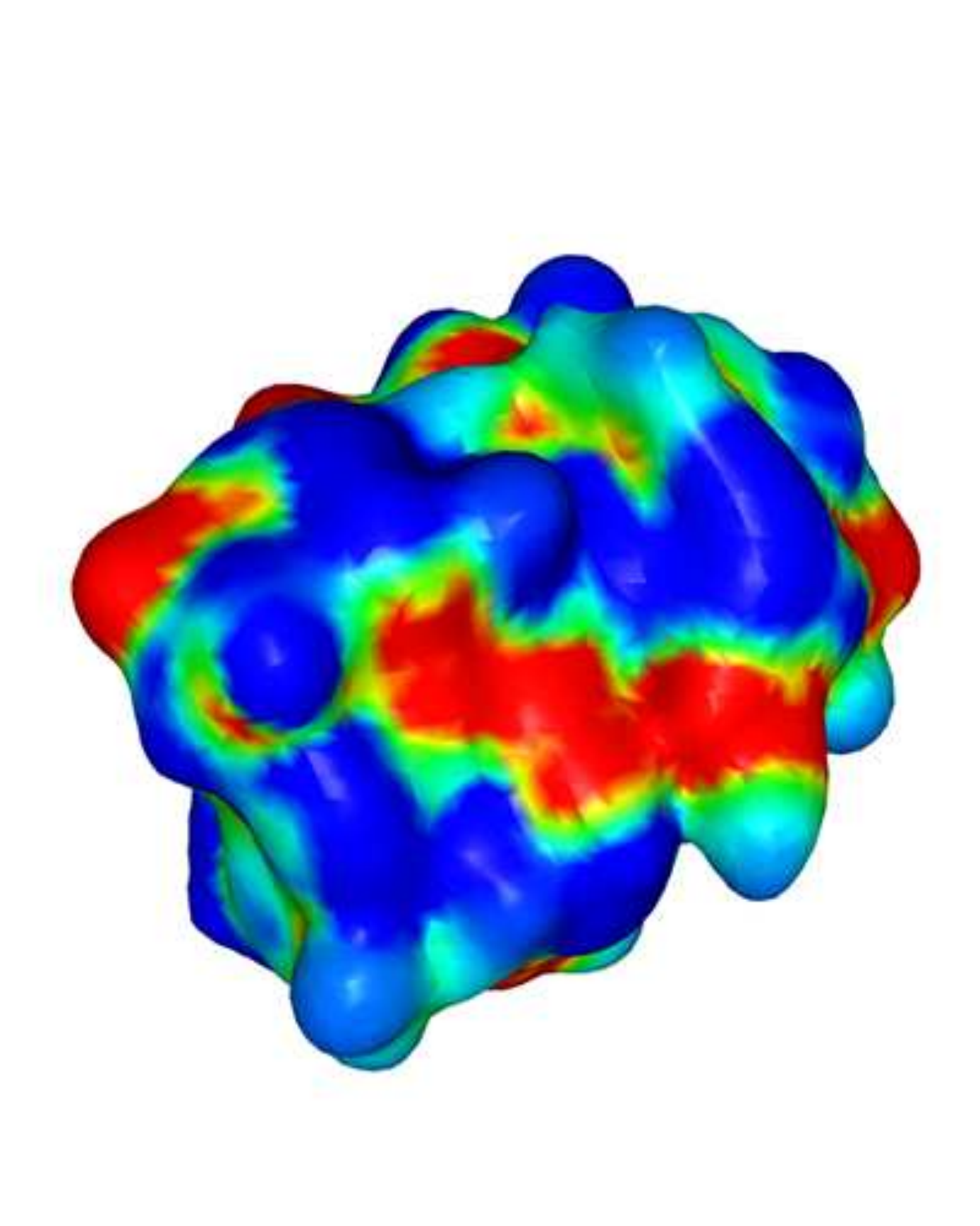}
  \includegraphics[scale = 0.18]{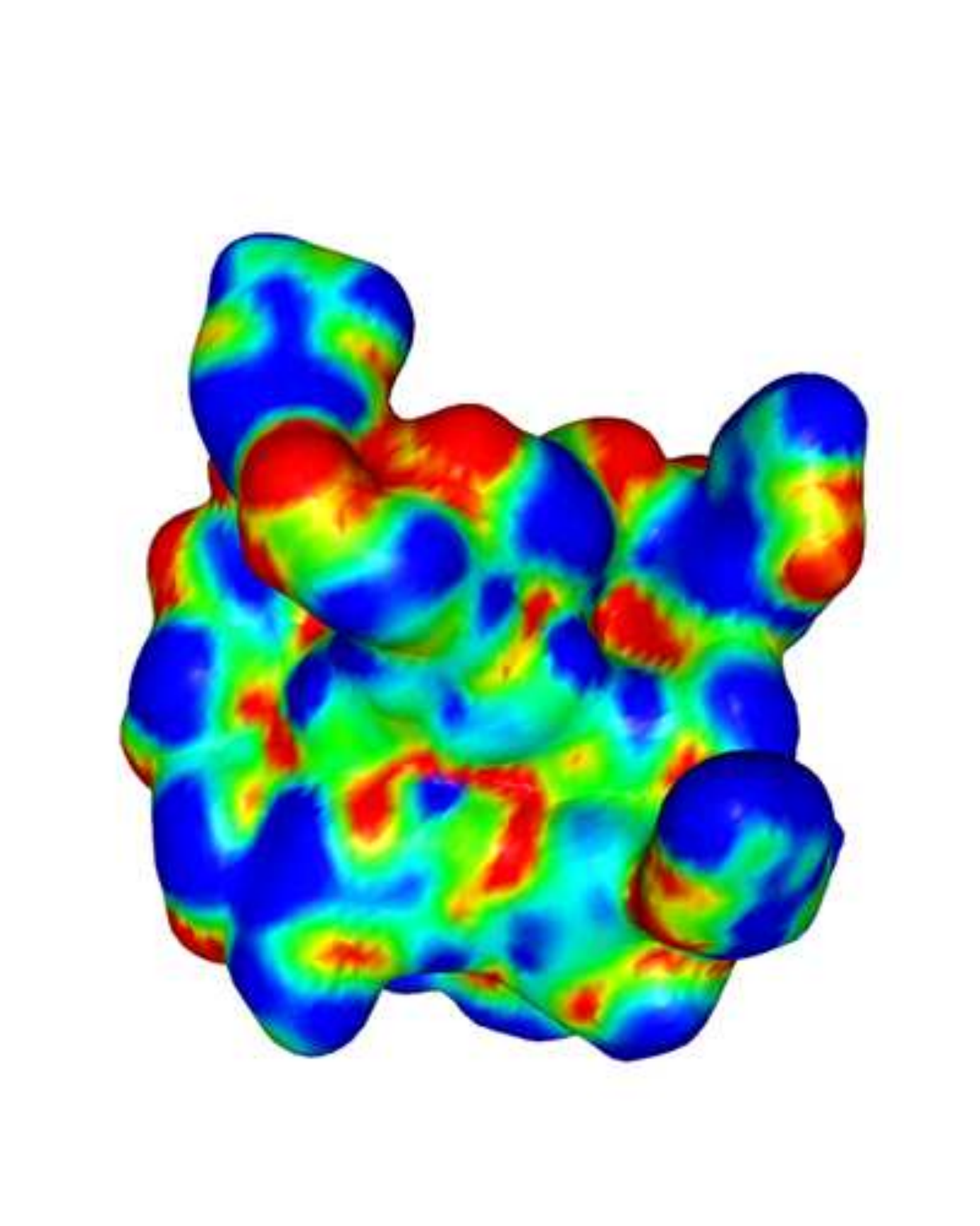}
  \includegraphics[scale = 0.18]{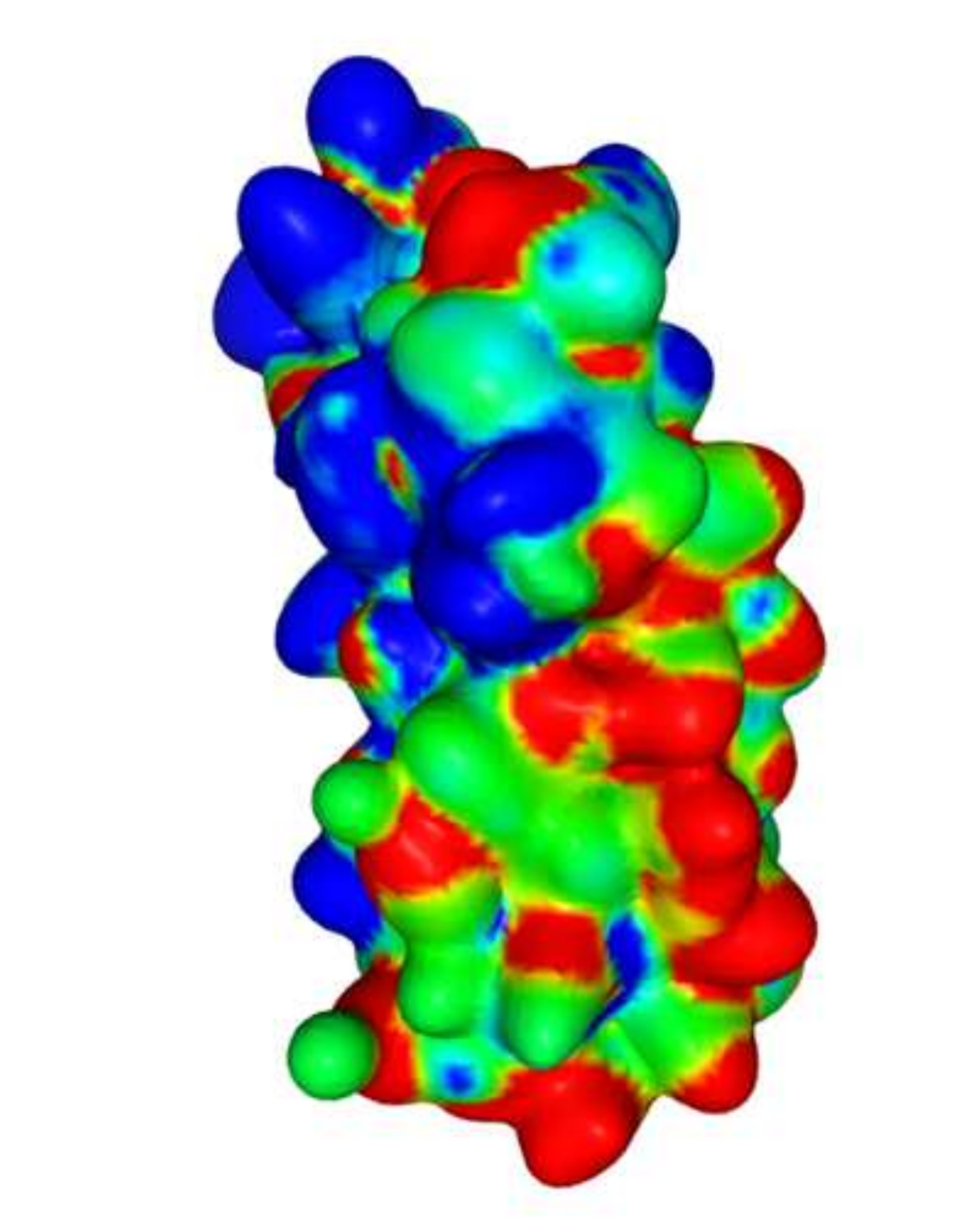}
  \caption{Electrostatic potential on molecular surfaces, calculated with AFMPB. From left to right: ADP, 2FLY, 6BST,  2O3M and 2IJI.}
  \label{pb}
\end{figure}

The algorithms introduced in the method section are used to generate the sparse surface. We here also test the applicability of the original and the sparse surface in the computations of Poisson-Boltzmann (PB) electrostatics. The boundary element method software used is a publicly available PB solver, AFMPB \cite{Zhang2015PAFMPB}. Table \ref{tab:solvationenergy} shows that AFMPB can undergo and produce converged results using the sparse represented surface, and the calculated solvation energies are close to the results using the original surface. Fig. \ref{pb}, using VISIM (www.xyzgate.com), shows the computed electro-static potentials mapped on the molecular surface. In the future, we can consider adding the charge information to the sparse representation.
\begin{table}[H]
\begin{center}
\caption{The solvation energy obtained with the original surface and the sparse represented surface for five biomolecules. Note: isovalue $\phi = 1.0$, initial decay rate $d = 0.5$.}\label{tab:solvationenergy}
\begin{tabular}{c|c|c|c}
\hline \hline
\multirow{2}*{Molecule}  & \multicolumn{3}{c}{Solvation energy $(kcal / mol)$}  \\
\cline{2-4}  
 & Original & Sparse & Relative error    \\  \hline
ADP  & -2.25992e+02  & -2.30075e+02  & 0.0181 \\ \hline
2FLY & -2.38927e+02  & -2.42670e+02   & 0.0157 \\ \hline
6BST & -9.16715e+02  & -9.20137e+02  & 0.0037 \\ \hline
2O3M & -3.03482e+03  & -3.05604e+03 & 0.0070 \\ \hline
2IJI &  -6.59502e+02  & -6.65894e+02  & 0.0097 \\ \hline \hline
\end{tabular}
\end{center}
\end{table}

\section{Conclusion}\label{section4}
In this paper, a sparse Gaussian molecular shape representation based on ellipsoid RBF neural network is proposed for arbitrary molecule. The original  Gaussian density maps is approximated with the ellipsoid RBF neural network. The sparsity of the ellipsoid RBF neural network is computed by solving an $L_1$ regularization optimization problem. Comparisons and experimental results indicate that our network needs much less number of ellipsoid RBF neurons to represent the original Gaussian density maps.

\section*{Acknowledgements}
This work was supported by the National Key Research and Development Program of China (grant 2016YFB0201304) and the China NSF (NSFC 11771435, NSFC 21573274).

\bibliographystyle{elsarticle-harv}
\bibliography{main}

%
%
%
\end{document}